\documentclass[12pt]{amsart}

\newcommand{\vsp}{\vspace{5mm}}
\newcommand{\bt}{\begin{Theorem}}
\newcommand{\et}{\end{Theorem}}
\newcommand{\bi}{\begin{itemize}}
\newcommand{\ei}{\end{itemize}}
\newcommand{\bea}{\begin{eqnarray}}
\newcommand{\eea}{\end{eqnarray}}
\newtheorem{Theorem}{\sc Theorem}
\newtheorem{Lemma}[Theorem]{\sc Lemma}
\newtheorem{Proposition}[Theorem]{\sc Proposition}
\newtheorem{Corollary}[Theorem]{\sc Corollary}
\newtheorem{Definition}[Theorem]{\sc Definition}

\newcommand{\be}{\begin{equation}}
\newcommand{\ee}{\end{equation}}

\def\qed{\hfill$\Box$}
\def\bi{\bibitem}

\def\CN{{\mathcal {N}}}
\def\CM{{\mathcal {M}}}
\def\CA{{\mathcal {A}}}
\def\CB{{\mathcal {B}}}
\def\CS{{\mathcal {S}}}
\def\CL{{\mathcal {L}}}
\def\CH{{\mathcal {H}}}
\def\CK{{\mathcal {K}}}

\def\CP{{\mathcal {P}}}
\def\CD{{\mathcal {D}}}

\def\CR{{\mathcal {R}}}
\def\CZ{{\rm\ke}rn.26em
\newcommand\la{{\langle}}
\newcommand\ra{{\rangle}}
\newcommand\lar{\leftarrow}
\newcommand\Lar{\Leftarrow}
\newcommand\rar{\rightarrow}
\newcommand\Rar{\Rightarrow}

\vrule width.02em height.5ex depth0ex \kern.04em \vrule width
.02em height1.47ex depth-1ex \kern-.34em Z}

\def\C{{\rm \kern.24em
 \vrule width.02em
    height1.4ex depth-.05ex
 \kern-.26em C}}

\def\ra{{\rightarrow}}

\def\ei{{\bf e_i}}

\def\uk{{\underline{k}}}
\def\ul{{\underline{l}}}

\def\uz{{\underline{z}}}
\def\ue{{\underline{e}}}
\def\oq{{\overline{q}}}
\def\uT{{\underline{T}}}
\def\uR{{\underline{R}}}
\def\uV{{\underline{V}}}

\def\uS{{\underline{S}}}

\def\circledq{{\bigcirc \!\!\!\!q}}

\def\N{{\rm I\kern-.23em N}}
\def\B{{\rm I\kern-.25em B}}
\def\D{{\rm I\kern-.25em D}}
\def\E{{\rm I\kern-.25em E}}
\def\F{{\rm I\kern-.25em F}}

\def\I{{\rm I\kern-.25em I}}

\def\M{{\rm I\kern-.23em M}}
\def\P{{\rm I\kern-.25em P}}
\def\A{{\rm \kern.22em
 \vrule width.02em
    height0.5ex depth 0ex
 \kern-.24em A}}
\def\G{{\rm \kern.24em
 \vrule width.02em
    height1.4ex depth-.05ex
 \kern-.26em G}}
\def\J{{\rm \kern.19em
 \vrule width.02em
    height1.47ex depth 0ex
 \kern-.21em J}}
\def\O{{\rm \kern.24em
 \vrule width.02em
    height1.4ex depth-0.5ex
 \kern-.26em O}}
\def\Q{{\rm \kern.24em
 \vrule width.02em
    height1.4ex depth-.05ex
 \kern-.26em Q}}
\def\S{{\rm \kern.18em
 \vrule width.02em
    height1.4ex depth-.9ex
  \kern.12em
  \vrule width.02em
     height0.7ex depth 0ex
  \kern-.34em S}}
\def\T{{\rm \kern.45em
 \vrule width.02em
    height1.47ex depth 0ex
 \kern-.47em T}}
\def\U{{\rm \kern.30em
 \vrule width.02em
    height1.47ex depth-.05ex
 \kern-.32em U}}
\def\V{{\rm \kern.27em
 \vrule width.02em
    height1.47ex depth-.8ex
 \kern-.29em V}}
\def\W{{\rm \kern.25em
 \vrule width.02em
    height1.47ex depth-0.9ex
 \kern.34em
 \vrule width.02em
    height1.47ex depth-.9ex
  \kern-.63em W}}
\def\X{{\rm \kern.30em
 \vrule width.02em
    height1.4ex depth-1ex
  \kern.12em
  \vrule width.02em
     height0.4ex depth 0ex
  \kern-.46em X}}
\def\Y{{\rm \kern.25em
 \vrule width.02em
    height1.0ex depth 0ex
 \kern-.27em Y}}
\def\Z{{\rm \kern.26em
 \vrule width.02em
    height0.5ex depth 0ex
  \kern.04em
  \vrule width.02em
     height1.47ex depth-1ex
  \kern-.34em Z}}

\begin{document}

\begin{center} {\bf {\Large Standard dilations
of $q$-commuting tuples}}

\end{center}

\vsp
\vsp
\vsp

\begin{center}

{\sc Santanu Dey } \vsp

{\bf October 24, 2003}
\end{center}

\vsp \vsp \vsp \vsp

\begin{center}
{\underline {\bf Abstract}}
\end{center}

\vsp

Here we study dilations of $q$-commuting tuples. In [BBD] the
authors gave the correspondence between the two standard dilations
of commuting tuples and here these results have been extended to
$q$-commuting tuples. We are able to do this when $q$-coefficients
`$q_{ij}$' are of modulus one. We introduce `maximal $q$-commuting
subspace ' of a $n$-tuple of operators and `standard q-commuting
dilation'. Our main result is that the maximal $q$-commuting
subspace of the standard noncommuting dilation of $q$-commuting
tuple is the `standard q-commuting dilation'. We also introduce
$q$-commuting Fock space as the maximal $q$-commuting subspace of
full Fock space and give a formula for projection operator onto
this space. This formula for projection helps us in working with
the completely positive maps arising in our study. The first
version of the Main Theorem (Theorem 19) of the paper for normal
tuples using some tricky norm estimates and then use it to prove
the general version of this theorem.

\vsp
\vfill
----------------------------------------------------------------------

\noindent {\sc Key words}: Dilation, $q$-Commuting Tuples,
Complete Positivity

\noindent {\sc Mathematics Subject Classification}: 47A45, 47A20

\newpage

\setcounter{equation}{0}
\begin{section}{Introduction}
A generalization of contraction operator in multivariate operator
theory is a contractive $n$-tuple which is defined as follows:
\begin{Definition}
{\em  A $n$-tuple $\underline{T}=(T_1,\ldots ,T_n)$ of bounded
operators on a Hilbert space ${\mathcal H}$ such that $T_1T_1^*+
\cdots  +T_nT_n^*\leq I$ is a {\em contractive $n$-tuple,\/} or a {\em row
contraction\/}.}
\end{Definition}

Along the lines of [BBD], we will study the dilation of a class of
operator tuples defined as follows:
\begin{Definition}
{\em A $n$-tuple $\uT =(T_1,\ldots, T_n)$ is said to be {\em $q$-commuting \/}if $T_j T_i =
q_{ij} T_i T_j$ for all $1 \leq i < j \leq n $, where $q_{ij}$ are
complex numbers.}
\end{Definition}
Such operator tuples appear often in Quantum Theory ([C] [M]
[Pr]). Here we introduce `maximal $q$-commuting piece' and using
this and a particular representation of permutation group we give
a definition for $q$-commuting Fock space when $q$-coefficients
`$q_{ij}$' are of modulus one. We have this condition for
$q$-coefficient for almost all the results here. This q-commuting
Fock space is different from the twisted Fock space of M.
Bo\.zejko and R. Speicher ([BS1]) or that of P. E.  T. Jorgensen
([JSW]). In section 2 we give formula for the projection of full
Fock space onto this space. We obtain a special tuple of
$q$-commuting operators and show that it is unitarily equivalent
to the tuple of shift operators of [BB]. We are able to show that
the range of the operator $A$ defined in equation (2.4) gives an
isometry onto the $q$-commuting Fock Space tensored with a Hilbert
space when $\uT$ is a  pure contractive tuple (this operator were
used by  Popescu and Arveson in [Po3], [Po4], [Ar2] and for
$q$-commuting case by Bhat and Bhattacharyya in [BB]). Using this
we are able to give a condition equivalent to the assertion of the
Main Theorem to hold for $q$-commuting purely contractive tuple.
In section 3 the proof of the particular case of Theorem 19 where
$\uT$ is also $q$-spherical unitary (introduced in section 3) is
more difficult than the version for commuting tuple and we had to
carefully choose the terms and proceed in a way that `$q_{ij}$' of
the $q$-commuting tuples get absorbed or cancel out when we
simplify the terms. Also unlike [BBD] we had to use an inequality
related to completely positive map before getting the result
through norm estimates. We are not able to generalize section 4 of
[BBD]. In the last section here we calculate the distribution of
$S_i+S_i^*$ with respect to the vacuum expectation and study some
properties of the related operator spaces.

For operator tuples $(T_1, \ldots, T_n)$, we need to consider the
products of the form $T_{\alpha _1}T_{\alpha _2}\cdots$ $
T_{\alpha _m}$, where each $\alpha _k\in \{ 1, 2, \ldots , n\}$.
We would have the following a notation for such products. Let
$\Lambda $ denote the set $\{ 1, 2, \ldots , n\}$ and $\Lambda ^m$
denote the $m$-fold cartesian product of $\Lambda $ for $m\geq 1.$
Given $\alpha =(\alpha _1, \ldots , \alpha _m)$ in $\Lambda ^m$,
$\uT ^{\alpha }$ will mean the operator $T_{\alpha _1}T_{\alpha
_2}\cdots T_{\alpha _m}$. Let $\tilde{\Lambda}$ denote
$\cup_{n=0}^{\infty} \Lambda^n$, where $\Lambda ^0$ is just the
set $\{ 0\}$ by convention and by $\uT ^0$ we would mean the
identity operator of the  Hilbert space where $T_i$'s are acting.

Let $\CS_m$ denote the group of permutation on $m$ symbols
$\{1,2,\cdots, m\}$. For a $q$-commuting tuple
$\underline{T}=(T_1,\ldots ,T_n),$ consider the product
$T_{x_1}T_{x_2}...T_{x_m}$ where $1 \leq x_i \leq n. $ If we
replace a consecutive pair say $T_{x_i}T_{x_{i+1}}$ of operators
in the above product by $q_{x_{i+1}x_i}T_{x_{i+1}}T_{x_i}$ and do
finite number of such operations with different choices of
consecutive pairs of these operators appearing in the subsequent
product of operators after each such operation, we will get a
permutation $\sigma \in \CS_m$ such that the final product of
operators can be written as
$kT_{x_{\sigma^{-1}(1)}}T_{x_{\sigma^{-1}(2)}}
...T_{x_{\sigma^{-1}(m)}}$ for some $k \in \C,$ i. \@ e.,
$T_{x_1}T_{x_2}...T_{x_m}=
kT_{x_{\sigma^{-1}(1)}}T_{x_{\sigma^{-1}(2)}}
...T_{x_{\sigma^{-1}(m)}}.$ For defining $q$-commuting tuple in
definition 2 we needed the known fact that this $k$ depends only
on $\sigma$ and $x_i,$ and not on the different choice of above
operations that give rise to the same final product of operators
$T_{x_{\sigma^{-1}(1)}}T_{x_{\sigma^{-1}(2)}}
...T_{x_{\sigma^{-1}(m)}}$. It also follows from the Proposition 5
in section 2.

\begin{Definition}
{\em Let ${\mathcal H}, {\mathcal L}$ be two Hilbert spaces such
that ${\mathcal H}$ be a closed subspace of ${\mathcal L}$ and let
 $\uT, \uR$ are $n$-tuples of bounded operators on
${\mathcal H}$, ${\mathcal L}$ respectively. Then $\uR $ is called
a {\em dilation\/} of  $\uT $ if
$$R_i^*u=T_i^*u $$
for all $u\in {\mathcal H}, 1\leq i\leq n.$ In such a case $\uT$
is called a {\em piece\/} of $\uR .$ If $\uT$ is a
$q$-commuting tuple ( i.\@ e., $T_jT_i=q_{ij} T_iT_j,$ for all
$i,j$), then it is called a $q$-commuting piece of $\uR .$ A
dilation $\uR $ of $\uT$ is said to be a {\em minimal dilation\/ }
if ${\overline {\mbox span}}\{ \uR ^{\alpha }h: \alpha \in {\tilde
\Lambda}, h\in \CH \}=\CL .$ And if $\uR$ is a tuple of
$n$ isometries and is a minimal dilation of $\uT$, then it is
called the {\em minimal isometric dilation \/} or the {\em
standard noncommuting dilation \/} of $\uT$.}

\end{Definition}

A presentation of the standard noncommuting dilation taken from [Po1]
has been used here to proof the main Theorem.
All Hilbert spaces that we consider will be complex and separable.
For a subspace $\CH$ of a Hilbert space, $P_{\CH}$ will denote the
orthogonal projection onto $\CH$. Standard noncommuting dilation
of $n$-tuple of bounded operators,
is unique upto unitary equivalence (refer [Po1-4]). Extensive
study of standard noncommuting dilation was carried out by Popescu.
He generalized many one variable results to multivariable case. It is easy to
see that if $\uR$ is a dilation of $\uT$ then
\begin{equation}
\uT ^{\alpha }(\uT ^{\beta })^*= P_{{\mathcal H}}\uR ^{\alpha
}(\uR ^{\beta })^*|_{{\mathcal H}},
\end{equation}
and for any polynomials $p, q$ in $n$-noncommuting variables
$$p(\uT)(q(\uT))^*=P_{{\mathcal H}}p(\uR )(q(\uR))^*|_{{\mathcal
H}}.$$ For a $n$-tuple $\uR $ of bounded operators on a Hilbert
space $\CM$, consider
 $${\mathcal C}^q (\uR )=\{ \CN : R_i ~~\mbox {leaves $\CN$ invariant },
~~R_i^*R_j^*h=\oq_{ij}R_j^*R_i^*h, \forall h\in \CN, \forall i,j\}.$$
 It is a complete lattice, in the sense that
arbitrary intersections and span closures of arbitrary unions of
such spaces are again in this collection. So it has a
maximal element and we denote it by $\CM ^q(\uR)$ (or by $\CM^q $
when the tuple under consideration is clear).

\begin{Definition}

{\em Let $\uR $ be a $n$-tuple of operators on a Hilbert space
$\CM$.  The $q$-commuting piece $\uR ^q=(R_1^q, \ldots , R_n^q)$
obtained by compressing $\uR$ to the maximal element $\CM ^q(\uR)$
of ${\mathcal C}^q (\uR)$ is called the {\em maximal $q$-commuting
piece\/ } of $\uR$. The maximal $q$-commuting piece is said to be
{\em trivial \/ } if  $\CM ^q(\uR)$ is the zero space.}
\end{Definition}

For any Hilbert space $\CK$, we have the full Fock space over
$\CK$ denoted by $\Gamma (\CK)$ as,
$$\Gamma (\CK)=\mathbb{C}\oplus \CK \oplus \CK^{\otimes ^2}\oplus \cdots
\oplus \CK ^{\otimes ^m}\oplus \cdots ,
$$
We denote the vacuum vector $1\oplus 0\oplus \cdots $  by
$\omega$.   For fixed
$n\geq 2,$ let ${\mathbb{C}}^n$ be the $n$-dimensional complex
Euclidian space with usual inner product and
$\Gamma({\mathbb{C}}^n)$ be the full Fock space over
${\mathbb{C}}^n$. Let $\{e_1, \ldots , e_n\}$ be the standard
orthonormal basis of ${\mathbb{C}}^n$. For $\alpha \in {\tilde \Lambda}$, $e^{\alpha }$ will
denote the vector $e_{\alpha _1}\otimes e_{\alpha _2}\otimes
\cdots \otimes e_{\alpha _m}$ in the full Fock space $\Gamma
(\mathbb{C}^n)$ and $e^0$ will denote the vacuum vector $\omega $. Then the (left) creation
operators $V_i$ on $\Gamma({\mathbb{C}^n})$ are defined by
\[ V_i x = e_i \otimes x \]  where $1 \leq i \leq n$ and $x
\in \Gamma({\mathbb{C}}^n)$ ( here $e_i\otimes \omega$ is
interpreted as $e_i$). It is obvious  that the tuple
$\underline{V}= (V_1, \ldots , V_n)$ consists of isometries with
orthogonal ranges and  $\sum V_iV_i^*=I-I_0$, where $I_0$ is the
projection on to the vacuum space. Let us define {\em
$q$-commuting Fock space} as the subspace $(\Gamma (\mathbb
{C}^n))^q (\uV)$ and let it be denoted by $\Gamma_q (\mathbb
{C}^n)$.

 Let $\uS =(S_1,
\ldots , S_n)$ be the tuple of operators on $\Gamma
_q(\mathbb{C}^n)$ where $S_i$ is the compression of $V_i$ to
$\Gamma _q(\mathbb{C}^n)$:
$$S_i=P_{\Gamma _q(\mathbb{C}^n)}V_i|_{\Gamma _q(\mathbb{C}^n)}.$$
Clearly each  $V_i^*$ leaves $\Gamma _q(\mathbb{C}^n)$ invariant.

Then it is easy to see that $\uS$ satisfies $\sum
S_iS_i^*=I^q-I^q_0$ (where $I^q, I^q_0$ are identity, projection
onto vacuum space respectively in $\Gamma _q(\mathbb{C}^n)$). So $\uV$ and
$\uS$ are contractive tuples, $S_jS_i=q_{ij}S_iS_j$ for all $1\leq
i, j\leq n,$ and $S_i^*x=V_i^*x$, for $x\in \Gamma
_q(\mathbb{C}^n)$.

The following result gives a description for maximal $q$-commuting
piece.

\begin{Proposition}
Let $\uR=(R_1,...,R_n)$ be a $n$-tuple of bounded operators on a
Hilbert space ${\mathcal M}$, $\CK _{ij}= \overline{span}\{ \uR
^{\alpha }(q_{ij} R_iR_j-R_jR_i)h: h\in \CM , \alpha \in \tilde
{\Lambda }\}$ for all $1\leq i, j\leq n$, and $\CK=\overline{span}
\{\cup_{i, j =1}^n \CK _{ij} \}$. Then $\CM ^q(\uR) = \CK ^{\perp
}$ and $\CM ^q(\uR)=\{ h\in \CM: (\oq_{ij}
R_j^*R_i^*-R_i^*R_j^*)(\uR ^\alpha )^*h=0, \forall 1\leq i, j\leq
n, \alpha\in {\tilde \Lambda}\}.$
\end{Proposition}
The above Proposition can be easily proved using arguement similar
to the proof of Proposition 4 of [BBD].

\begin{Corollary}
Suppose $\uR $, $\uT $ are $n$-tuples of operators on two Hilbert
spaces $\CL , \CM$. Then the maximal $q$-commuting piece of
$(R_1\oplus T_1, \ldots , R_n\oplus T_n)$ acting on $\CL\oplus
\CM$ is $(R_1^q \oplus T_1^q, \ldots , R_n^q\oplus T_n^q)$ acting
on $\CL ^q\oplus \CM ^q$ and the maximal $q$-commuting piece of
$(R_1\otimes I, \ldots , R_n\otimes I)$ acting on $\CL \otimes
\CM$ is $(R_1^q\otimes I, \ldots , R_n^q\otimes I)$ acting on $\CL
^q\otimes \CM.$
\end{Corollary}

\noindent {\sc Proof:} Clear from Proposition 6. \qed

\begin{Proposition}
Let $\uT, \uR$ are $n$-tuples of bounded operators on ${\mathcal
H}$, ${\mathcal L}$, with $\CH \subseteq \CL$, such that $\uR$ is
a dilation of $\uT$. Then $\CH ^q(\uT)=\CL ^q(\uR)\bigcap \CH$ and
$\uR ^q$ is a dilation of $\uT ^q$.
\end{Proposition}

\noindent {\sc Proof:} This can be using arguements similar to
proof of Proposition 7 of [BBD].  \qed
\end{section}

\begin{section}{A $q$-Commuting Fock Space}

\setcounter{equation}{0}

For a $q$-commuting $n$-tuple $\uT$ on a finite dimensional Hilbert space $\CH$
say of dimension $m,$ because of the relation
$$\mbox{Spectrum} (T_iT_j)\cup \{0\}=\mbox{Spectrum} (T_jT_i)
\cup \{0\}=\mbox{Spectrum} (q_{ij}T_iT_j)\cup \{0\},$$
we get $q_{ij}$ is either $0$ or $m^{\mbox{th}}$-root of unity.

Here after whenever we deal with $q$-commuting tuples we would
have another condition on the tuples that $|q_{ij}|=1$ for $1\leq
i, j \leq n.$ However Proposition 5, Proposition 6 and Corollary 7 does not
need this assumption. Let $\underline{T}=(T_1,\ldots ,T_n)$ be a $q$-commuting tuple
and consider the product $T_{x_1}T_{x_2}...T_{x_m}$ where $1 \leq
x_i \leq n. $ Let  $\sigma \in \CS_m.$ As transpositions of the type
$(k,k+1), 1 \leq k \leq m-1$ generates $\CS_m,$  Let $\sigma^{-1}$ be $ \tau_1 \ldots
\tau_s$ where for each $1 \leq i \leq s$ there exist
 $k_i$ such
that $1 \leq k_i \leq m-1$ and $\tau_i$ is a transposition of the
form $(k_i, k_i+1).$  Let $\tilde{\sigma}_i = \tau_{i+1} \tau_{i}
\ldots \tau_{s}$ for $1 \leq i \leq s-1$ and $\tilde{\sigma}_s$ be
the identity permutation. Let us define $y_i=x_{\tilde{\sigma}_{i}(k_{i})}$
and $z_i=x_{\tilde{\sigma}_{i}(k_{i}+1)}$. If we substitute
$T_{y_s}T_{z_s}$ by
$q_{z_sy_s}T_{z_s}T_{y_s}$ corresponding to
$\tau_s,$ substitute
$T_{y_{s-1}}T_{z_{s-1}}$ by $q_{z_{s-1}y_{s-1}} T_{z_{s-1}}T_{y_{s-1}}$
corresponding to $\tau_{s-1},$
 and so on till we substitute the corresponding
term for $\tau_1,$ we would get
$q^\sigma_1(x) \ldots
q^\sigma_s(x)T_{x_{\sigma^{-1}(1)}}
T_{x_{\sigma^{-1}(2)}} \cdots T_{x_{\sigma^{-1}(m)}}$  where
$q^{\sigma}_i(x) =
q_{z_iy_i}.$ That is
$T_{x_1}T_{x_2} \cdots  T_{x_m}=q^\sigma_1(x) \ldots
q^\sigma_s(x)T_{x_{\sigma^{-1}(1)}}
T_{x_{\sigma^{-1}(2)}} \cdots  T_{x_{\sigma^{-1}(m)}}.$ Let
$q^{\sigma}(x) = q^\sigma_1(x) \ldots q^\sigma_s(x)$
where $q^{\sigma}_i(x) = q_{z_iy_i}$.

\begin{Proposition}
Let $\underline{T}=(T_1,\ldots ,T_n)$ be a $q$-commuting tuple and
consider the product $T_{x_1}T_{x_2}...T_{x_m}$ where $1 \leq x_i
\leq n. $ Suppose  $\sigma \in \CS_m$ and $ q^\sigma (x)$ be as
defined above. Then
$$q^{\sigma}(x)=\prod q_{x_{\sigma^{-1}(k)}x_{\sigma^{-1}(i)}},$$
where product is over $\{(i,k): 1\leq i < k \leq m, \sigma^{-1}(i)
> \sigma^{-1}(k) \}.$ Moreover $q^{\sigma}(x)$ does not depend on
the choice of $\sigma.$
\end{Proposition}

\noindent {\sc Proof:}  We have
$$q^{\sigma}(x) = q^\sigma_1(x) \ldots q^\sigma_s(x)$$
where $q^{\sigma}_i(x) = q_{z_iy_i}$.
 For a pair $i,k$ such that $1\leq i < k \leq
m$ let $k'=\sigma^{-1}(k)$ and $i'=\sigma^{-1}(i).$ Let $\sigma=
\tau_1 \cdots, \tau_s$ and $\tilde{\sigma}_i$ be as defined above.
If $i'
> k'$ then there are odd number of transpositions $\tau_r$ for $1
\leq r \leq m$ such that they interchange the positions of $i'$
and $k'$ in the image of $\tilde{\sigma}_r$ when we consider the
composition $\tau_r \tilde{\sigma}_r.$ And for $1\leq i < k \leq
m$ if $i' < k'$ then there are even number of transpositions
$\tau_r$ for $1 \leq r \leq m$ such that they interchange the
positions of $i'$ and $k'$ in the image of $\tilde{\sigma}_r$ when
we consider the composition $\tau_r\tilde{\sigma}_r.$  For the
first transposition in $\tau_r$ that interchanges $i'$ and $k',$
the corresponding factor in $q^\sigma(x)$ say $q^\sigma_r(x)$ is
$q_{x_{k'}x_{i'}},$ for the second transposition that interchanges
$i'$ and $k',$
 the corresponding factor is
$q_{x_{i'}x_{k'}},$ for the third transposition that interchanges $i'$ and $k',$
the corresponding factor is
$q_{x_{k'}x_{i'}},$ and so on. But
$(q_{x_{i'}x_{k'}})^{-1}=q_{x_{k'}x_{i'}}$ and so
$$q^\sigma(x)= \prod q_{x_{\sigma^{-1}(i)}x_{\sigma^{-1}(k)}},$$ where
product is over $\{(i,k): 1\leq i < k \leq m, \sigma^{-1}(i) >
\sigma^{-1}(k) \}.$\qed

Following similar arguements it is easy to see
that if there exist $\sigma \in \CS_m$ such that $(x_1,\cdots,x_n)
=(x_{\sigma^{-1}(1)},\cdots,x_{\sigma^{-1}(n)}),$ then $q^\sigma (x)=1.$

 Let
$U^{m,q}_{\sigma}$ be defined on $(\mathbb {C}^n)^{\otimes ^m}$ by

\begin{equation}
U^{m,q}_{\sigma}( e_{x_1} \otimes \ldots \otimes e_{x_m} ) =
q^{\sigma}(x)e_{x_{\sigma^{-1}(1)}} \otimes \ldots \otimes
e_{x_{\sigma^{-1}(m)}}
\end{equation}

 on the standard basis vectors and extended
linearly on $(\mathbb {C}^n)^{\otimes ^m}$. As $|q_{ij}|=1$ for $1
\leq i,j \leq n,~ $ $U^m_{\sigma}$ is unitary and $U^m_{\sigma}$
extends uniquely to a unitary operator on $(\mathbb
{C}^n)^{\otimes ^m}$.

 Let $$(\mathbb{C}^n)^{\circledq ^m} = \{u \in (\mathbb
{C}^n)^{\otimes ^m}: U^{m,q}_{\sigma} u = u ~ \forall \sigma \in
\CS_m \}$$ and $(\mathbb {C}^n)^{\circledq ^0} = \mathbb{C}$

\begin{Lemma}
 The map defined from $\CS _m$ to $B(({\mathbb C}^n)^{\otimes^m})$
defined by $\sigma \to U^{m,q}_{\sigma}$ for all $\sigma \in \CS_m
$ is a representation.
\end{Lemma}

 \noindent {\sc Proof:} Let $\otimes_{i=1}^m
e_{x_i}, \otimes_{i=1}^m e_{y_i} \in (\mathbb {C}^n)^{\otimes
^m}, 1 \leq x_i, y_i \leq n.$ Suppose there exist $\sigma \in \CS_m$ such that
$\otimes_{i=1}^m e_{y_i} = \otimes_{i=1}^m
e_{x_{\sigma^{-1}(i)}}.$ Then
 $\langle U^{m,q}_{\sigma}(\otimes_{i=1}^m e_{x_i}), \otimes_{i=1}^m e_{y_i}
\rangle = q^\sigma(x)$  and  $\langle \otimes_{i=1}^m e_{x_i},
U^{m,q}_{\sigma^{-1}}(\otimes_{i=1}^m
e_{y_i})\rangle=\overline{q^{({\sigma}^{-1})}(y)}.$ Also
$$q^{(\sigma^{-1})}(y)=\prod q_{y_{\sigma(k)}y_{\sigma(i)}}
= \prod q_{x_{k}x_{i}} $$ where the products are over $\{(i,k):
1\leq i < k \leq m,  \sigma(i) > \sigma(k) \}.$ If we substitute
$k=\sigma^{-1}(i')$  and  $i=\sigma^{-1}(k')$ in the last term we
get
$$q^{(\sigma^{-1})}(y)= \prod q_{x_{\sigma^{-1}(i')}x_{\sigma^{-1}(k')}}=(\prod
q_{x_{\sigma^{-1}(k')}x_{\sigma^{-1}(i')}})^{-1}=(q^\sigma (x))^{-1}$$
where the products are
over $\{(i',k'): 1\leq i'< k' \leq m, \sigma^{-1}(i')>
\sigma^{-1}(k') \}.$ So
$$q^\sigma(x)=(q^{({\sigma}^{-1})}(y))^{-1}=\overline{q^{({\sigma}^{-1})}(y)}.$$
The last equality holds as $|q_{ij}|=1.$ This implies $\langle
U^{m,q}_{\sigma}(\otimes_{i=1}^m e_{x_i}), \otimes_{i=1}^m e_{y_i}
\rangle = \langle \otimes_{i=1}^m e_{x_i},
U^{m,q}_{\sigma^{-1}}(\otimes_{i=1}^m e_{y_i})\rangle.$ If there
does not exist any $\sigma \in \CS_m$ such that $\otimes_{i=1}^m
e_{y_i} = \otimes_{i=1}^m e_{x_{\sigma^{-1}(i)}}$ then
$$\langle U^{m,q}_{\sigma'}(\otimes_{i=1}^m e_{x_i}), \otimes_{i=1}^m
e_{y_i} \rangle = 0 = \langle \otimes_{i=1}^m e_{x_i},
U^{m,q}_{(\sigma')^{-1}}(\otimes_{i=1}^m e_{y_i})\rangle$$ for all
$\sigma' \in \CS_m.$ So
$(U^{m,q}_{\sigma})^*=U^{m,q}_{\sigma^{-1}}$ for $\sigma \in
\CS_m,$ when acting on the basis elements of the
$(\mathbb{C}^n)^{\otimes^m},$ and hence is true for all elements
$(\mathbb{C}^n)^{\otimes^m}.$

Next let  $\sigma \in \CS_m$ be equal to $\sigma_1 \sigma_2$ for
some  $\sigma_1, \sigma_2 \in \CS_m.$ We would show that
$U^{m,q}_{\sigma}=U^{m,q}_{\sigma_1}U^m_{\sigma_2}.$ Let $e_x=
e_{x_1} \otimes \ldots \otimes e_{x_m}$ where $x_j \in
\{1,...,n\}$ for $1 \leq j \leq m.$ Let $\sigma_1^{-1} = \tau_1
\ldots \tau_r$ and $\sigma_2^{-1} = \tau_{r+1} \ldots \tau_s$
where for each $1 \leq i \leq s$, there exist $k_i$ such that $1
\leq k_i \leq m-1$ and $\tau_i$ is a transposition of the form
$(k_i, k_i+1).$
$$ U^{m,q}_{\sigma_1}U^{m,q}_{\sigma_2}( e_{x_1} \otimes \ldots \otimes e_{x_m} )
=U^{m,q}_{\sigma_1} (q^{\sigma_2}(x)e_{x_{\sigma_2^{-1}(1)}}
\otimes \ldots \otimes e_{x_{\sigma_2^{-1}(m)}}) = q^{\sigma_1}(z)
q^{\sigma_2}(x)e_{x_{\sigma_2^{-1}\sigma_1^{-1}(1)}} \otimes
\ldots \otimes e_{x_{\sigma_2^{-1}\sigma_1^{-1}(m)}}$$ where
$e_z=e_{z_1} \otimes \ldots \otimes e_{z_m},$ i.\@e,
$z_i=x_{\sigma_2^{-1}(i)}.$ But as $\sigma=\tau_1 \ldots \tau_r
\tau_{r+1} \ldots \tau_s$ it is easy to see that
$q^{\sigma}(x)=q^{\sigma_1}(z)q^{\sigma_2}(x)$ using the
definition of $q^{\sigma}(x).$ So we get
$$ U^{m,q}_{\sigma_1}U^{m,q}_{\sigma_2}( e_{x_1} \otimes \ldots \otimes e_{x_m} )
= q^{\sigma}(x) e_{x_{\sigma^{-1}(1)}} \otimes \ldots \otimes
e_{x_{\sigma^{-1}(m)}} = U^{m,q}_{\sigma}( e_{x_1} \otimes \ldots
\otimes e_{x_m} ). $$ And hence $U^{m,q}_{\sigma_1
\sigma_2}=U^{m,q}_{\sigma_1} U^{m.q}_{\sigma_2}.$
  \qed

Now if we use $Y^{m,q}_{\sigma}$ also to denote a operator in
$\Gamma ({\mathbb C}^n)$ which acts as $U^{m,q}_{\sigma}$ on
$(\mathbb {C}^n)^{\otimes ^m}$ and $I$ on the orthogonal, we get a
representation of $S_m$ on $B(\Gamma(\mathbb{C}^n)).$ In the next
Lemma and Proposition we derive a formula for the projection
operator onto the $q$-commuting Fock space.

\begin{Lemma}
Let $P_m$ be a operator on $(\mathbb {C}^n)^{\otimes ^m}$
defined by
\begin{equation}
P_m = \frac{1}{m !} \sum_{\sigma \in \CS_m} U^{m,q}_{\sigma}.
\end{equation}
Then $P_m $ is a projection of $(\mathbb
{C}^n)^{\otimes ^m}$ onto $(\mathbb {C}^n)^{\circledq ^m}$.
\end{Lemma}

\noindent {\sc Proof:}  First we see that
$$ P_m^* = \frac{1}{m !} \sum_{\sigma \in \CS_m}
(U^{m,q}_{\sigma})^* = \frac{1}{m !} \sum_{\sigma \in \CS_m}
U^{m,q}_{\sigma^{-1}} = P_m,$$
 Consider a permutation  $\sigma' \in \CS_m.$
$$ P_m U^{m,q}_{\sigma'} =  \frac{1}{m !} \sum_{\sigma \in \CS_m}
U^{m,q}_{\sigma\sigma'}= \frac{1}{m !} \sum_{\sigma \in \CS_m}
U^{m,q}_{\sigma}= P_m.$$ Similarly $ U^{m,q}_{\sigma'} P_m = P_m.$
So $P^2_m=P_m$ and hence $P_m$ is a projection. \qed

\begin{Proposition}
$\oplus_{m=0}^{\infty} (\mathbb {C}^n) ^{\circledq ^m} = \Gamma_q
(\mathbb {C}^n)$
\end{Proposition}

\noindent {\sc Proof:} Let $Q=\oplus_{m=0}^{\infty} P_{m}$ be the
projection of $\Gamma (\mathbb {C}^n)$ onto $\oplus_{m=0}^{\infty}
(\mathbb {C}^n) ^{\circledq ^m}$ where $P_m$ is a defined in Lemma
9. For transposition $(1,2),$ let us define $U^q_{(1,2)}$ as
$\oplus_{m=o}^{\infty} U^{m,q}_{(1,2)}$ where $U^{0,q}_{(1,2)}=I$
and $U^{1'q}_{(1,2)}=I.$ Let $\otimes_{i=1}^k e_{x_i}  \in
(\mathbb {C}^n)^{\otimes ^k}, 1 \leq x_i \leq n.$ Then
$$U^q_{(1,2)}
 V_jV_i (\otimes_{i=1}^k e_{x_i})= U^q_{(1,2)} \{e_j \otimes e_i \otimes
 (\otimes_{i=1}^k e_{x_i})\}  =q_{ij} e_i \otimes e_j \otimes
 (\otimes_{i=1}^k e_{x_i})=q_{ij} V_iV_j (\otimes_{i=1}^k e_{x_i}).$$

Next we would show that
 $\oplus_{m=0}^{\infty} (\mathbb {C}^n) ^{\circledq ^m}$ is
left invariant by $V_i^*.$ Let  $\otimes_{j=1}^m
e_{x_j} \in (\mathbb {C}^n)^{\otimes ^m},1 \leq x_j \leq n.$
 Then $V_i^*\{P_m (\otimes_{j=1}^m e_{x_j})\}$
is zero if none of $x_j$ is equal to $i.$ Otherwise $V_i^* \{ P_m ( \otimes_{j=1}^m e_{x_j})\}$
 is some non-zero element belonging to  $\oplus_{m=0}^{\infty} (\mathbb {C}^n)
^{\circledq ^{(m-1)}}$ because of the following.
Let  $x_j=i$ iff $j \in \{i_1,...,i_p\},$ and let $\CA_k$ be the set of all
$\sigma \in S_m$ such that $\sigma^{-1}$ sends $1$ to $i_k, 1 \leq k \leq p,$ then each
element of $\CA_k$ is a composition $\tau \sigma'$ where $\tau$ is the transposition $(1,i_k)$
and a permutation $\sigma'$ for which $({\sigma'})^{-1}$ keeps $1$ fixed and permutes
rest of the $m-1$ symbols. As $V_i$ are isometries with orthogonal ranges,
\begin{eqnarray*}
V_i^*\{P_m (\otimes_{j=1}^m e_{x_j})\} &=&V_i^*\{\frac{1}{m !}
\sum_{\sigma \in \CS_m} U^{m,q}_{\sigma}(\otimes_{j=1}^m
e_{x_j})\}
=\frac{1}{m !} \sum_{k=1}^p V_i^*(\sum_{\sigma \in \CA_{i_k}} U^{m,q}_{\sigma}e_{x_j})\\
&=&\frac{1}{m !} \sum_{k=1}^p V_i^*\{\sum_{\tau \sigma' \in \CA_{i_k}}
U^{m,q}_{\tau}U^{m,q}_{\sigma'} (\otimes_{j=1}^m e_{x_j})\} \\
&=& \sum_{k=1}^p  a_k (x) P_{m-1} (\otimes_{j=1}^m
e_{x_1} \otimes \cdots \otimes  \hat{e}_{x_{i_k}} \otimes \cdots \otimes e_{x_{m}})
\end{eqnarray*}
where $a_k(x)$ are constants and $\hat{e}_{x_{p}}$ denotes the term
 $e_{x_1} \otimes \cdots \otimes e_{x{p-1}} \otimes e_{x_{p+1}} \otimes \cdots \otimes e_{x_{m}}.$ This shows that
 $\oplus_{m=0}^{\infty} (\mathbb {C}^n) ^{\circledq ^m}$ is
left invariant by $V_i^*.$

Using these and the results of Lemma 9 we have the following.
Taking $R_i=QV_iQ$ for $\alpha \in \Lambda^m$ we get
\begin{eqnarray*}
 R_i R_j R^\alpha \omega = QV_iV_j V^\alpha \omega = QU^{m+2,q}_{(1,2)}
 q_{ji}V_jV_i V^\alpha \omega = q_{ji} QV_jV_i V^\alpha \omega
 = q_{ji} (R_j)(R_i) R^\alpha \omega.
\end{eqnarray*}
So $( QV_1Q, \ldots, QV_n Q )$ is a $q$-commuting piece of $\uV$.
To show maximality we make use of Proposition 6. Suppose $x \in
\Gamma ({\mathbb{C}}^n)$ and $\langle x, \uV
^\alpha(q_{ij}V_iV_j-V_jV_i)y\rangle =0$ for all $\alpha \in
\tilde{\Lambda }, 1\leq i,j\leq n$ and $y\in \Gamma
({\mathbb{C}}^n)$. We wish to show that $x\in \Gamma
_q({\mathbb{C}}^n)$. Suppose $x_m$ is the $m$-particle component
of $x$, i.\@e., $x=\oplus _{m\geq 0}x_m$ with $x_m \in
(\mathbb{C}^n)^{{\otimes}^m}$ for $m\geq 0$. For $m\geq 2$ and any
permutation $\sigma $ of $\{1,2, \ldots ,m\}$ we need to show that
the unitary $U^{m,q}_{\sigma }:(\mathbb{C}^n)^{{\otimes}^m} \to
(\mathbb{C}^n)^{{\otimes}^m}$, defined by equation (2.1) leaves
$x_m$ fixed. Since $\CS_m$ is generated by the set of
transpositions $\{(1,2), \ldots , (m-1, m)\}$ it is enough to
verify $U^{m,q}_{\sigma }(x_m)=x_m$ for permutations $\sigma $ of
the form $(i, i+1)$. So fix $m$ and $i$ with $m\geq 2$ and $1\leq
i\leq (m-1).$  We have
\begin{equation}
  \langle \oplus _px_p, \uV ^\alpha(q_{kl}V_kV_l-V_lV_k)\uV^\beta
\omega  \rangle =0,
\end{equation}
for every $\beta \in \tilde{\Lambda}, 1\leq k,l\leq n.$ This implies that
$$\langle x_m, e^\alpha \otimes (q_{kl}e_k\otimes e_l- e_l\otimes e_k)
\otimes e^\beta \rangle =0$$
for any $\alpha \in \Lambda^{i-1}, \beta \in \Lambda^{m-i-1}.$ So if
$$x_m=\sum a(s,t,\alpha,\beta )e^\alpha
\otimes e_s \otimes e_t \otimes e^\beta$$ where the sum is over
$\alpha \in \Lambda^{i-1}, \beta \in \Lambda^{m-i-1}$ and $1\leq
s,t\leq n,$ and $a(s,t,\alpha,\beta )$ are constants,then for
fixed $\alpha$ and $\beta$ it follows from equation (2.3) that
$\oq_{kl}a(k,l,\alpha,\beta )=a(l,k,\alpha,\beta )$ or
$q_{lk}a(k,l,\alpha,\beta )=a(l,k,\alpha,\beta ).$ Therefore
\begin{eqnarray*}
U^{m,q}_\sigma ( a(k,l,\alpha,\beta )e^\alpha \otimes e_k \otimes
e_l \otimes e^\beta +a(l,k,\alpha,\beta )e^\alpha
\otimes e_l \otimes e_k \otimes e^\beta ) \\
 = q_{lk}a(k,l,\alpha,\beta )e^\alpha \otimes e_l \otimes e_k \otimes e^\beta
+q_{kl}a(l,k,\alpha,\beta )e^\alpha \otimes e_k \otimes e_l \otimes e^\beta\\
 = a(l,k,\alpha,\beta )e^\alpha
\otimes e_l \otimes e_k \otimes e^\beta+ a(k,l,\alpha,\beta )e^\alpha \otimes e_k \otimes e_l \otimes e^\beta
\end{eqnarray*}

This clearly implies $U{m,q}_{\sigma }(x_m)=x_m$, for $\sigma =(i,
i+1)$. \qed

Let $\CP$ be the vector space of all polynomials in $q$-commuting
variables $z_1, \ldots, z_n$ that is $z_jz_i=q_{ij}z_iz_j.$ Any
multi-index $\uk$ is a ordered $n$-tuple of non-negative integers
$(k_1, \ldots,k_n)$. We shall write $k_1+ \ldots + k_n$ as
$|\uk|$. The special multi-index which has $0$ in all positions
except the $i^{th}$ one, where it has $1$, is denoted by $\ue_i$.
For any non-zero multi-index $\uk$ the monomial $z_1^{k_1} \ldots
z_n^{k_n}$ will be denoted by $\uz^{\uk}$ and for the multi-index
$\uk= (0, \ldots, 0),$ let $\uz^{\uk}$ be the complex number $1$.
Let us have the following inner product with it. Declare
$\uz^{\uk}$ and $\uz^{\ul}$ orthogonal if $\uk$ is not the same as
$\ul$ as ordered multi-indices. Let $$\| \uz^{\uk} \|^2 =
\frac{k_1 ! \cdots k_n ! }{|k|!}.$$ Note that the following
inner-product is also refered in [BB] in Definition (1.1) in
general case. Now define $\CH'$ to be the closure of $\CP$ with
respect to this inner product. Define a tuple $\uS'= (S_1',
\ldots, S_n')$ where each $S_i'$ is defined for $f \in \CP$ by
$$ S_i' f(z_1, \ldots, z_n) = z_i f(z_1, \ldots, z_n)$$
and $S_i$ is linearly extended to $\CH'$. In the case of our standard
$q$-commuting $n$-tuple $\uS$ of operators on $\Gamma_q (\mathbb
{C}^n)$, when $\uk=(k_1, \ldots,k_n)$  let $\uS^{\uk}= S_1^{k_1}
\ldots S_n^{k_n}$ and when $\uk= (0, \ldots, 0)$ let $\uS^{\uk} = 1.$

Using (2.2) and  the fact that  $V_i$'s are isometries
with orthogonal ranges for $\uk=(k_1, \ldots,k_n),|\uk|=m$  we get
 $$\|\uS^{\uk}\omega\|= \langle P_m \uV^{\uk} \omega, \uV^{\uk} \omega \rangle =
\langle \frac{1}{|\uk| !} \sum_{\sigma \in \CS_m} U^{m,q}_{\sigma}
\uV^{\uk} \omega, \uV^{\uk} \omega \rangle = \frac{k_1 ! \cdots
k_n !}{|\uk|!}.$$ If we denote $\uV^{\uk} \omega$ by $e_{x_1}
\otimes \cdots \otimes e_{x_m}, 1 \leq x_i \leq n,$ then to get
the last term of the above equation we used the fact that there
are $k_1 ! \cdots k_n !$ permutations $\sigma \in \CS_m$ such that
$e_{x_1} \otimes \cdots \otimes e_{x_m} = e_{x_{\sigma^{-1}(1)}}
\otimes \cdots \otimes e_{x_{\sigma^{-1}(m)}}$ .
  Next we show that the above tuples $\uS'$
and $\uS$ are unitarily equivalent.

\begin{Proposition} Let $\uS'= (S_1', \ldots, S_n')$
be the operator tuples on $\CH'$ as introduced above and let $\uS=
(S_1, \ldots, S_n)$ be the  standard $q$-commuting tuple of
operators on $\Gamma_q (\mathbb {C}^n)$. Then there exist unitary
$U: \CH' \to \CH$ such that $US_i'= S_i U$ for $1 \leq i \leq n$.
\end{Proposition}
\noindent {\sc Proof :} Define $U: {\CP} \to \Gamma_q (\mathbb
{C}^n)$ as
$$ U ( \sum_{|\uk| \leq s } b_{\uk} \uz^{\uk}) = \sum_{ |\uk| \leq
s} b_{\uk} \uS^{\uk}\omega $$
where $b_{\uk} \uz^{\uk} \in \CP,$ $b_{\uk}$ are constants.
As $\|\uz^{\uk} \|=\| \uS^{\uk}\omega \|$
we have $$\|\sum_{|\uk| \leq s } b_{\uk} \uz^{\uk} \|^2= \sum_{|\uk| \leq s }
|b_{\uk}|^2 \|\uz^{\uk} \|^2 = \sum_{|\uk| \leq s } |b_{\uk}|^2 \| \uS^{\uk}\omega \| ^2
= \|\sum_{|\uk| \leq s } b_{\uk} \uS^{\uk}\omega \|^2.$$
So we can extend it
linearly to ${\mathcal H}'$ and $U$ is a unitary.
\begin{eqnarray*}
U S_i'(\sum_{|\uk| \leq s } b_{\uk} \uz^{\uk})
&=& U (z_i\sum_{|\uk| \leq s } b_{\uk} \uz^{\uk})= q_{1i}^{k_1} \cdots q_{i-1i}^{k_{i-1}}
U(\sum_{|\uk| \leq s } b_{\uk} \uz^{\uk + \ue_i})\\
&=& q_{1i}^{k_1} \cdots q_{i-1i}^{k_{i-1}} \sum_{|\uk| \leq s } b_{\uk} \uS^{\uk + \ue_i}
\omega = S_i (\sum_{|\uk| \leq s } b_{\uk} \uS^{\uk} \omega )\\
&=& S_i U (\sum_{|\uk| \leq s } b_{\uk} \uz^{\uk}),
\end{eqnarray*}
i.\@ e., $U S_i'= S_i U$  for $1 \leq i \leq n$.
 \qed

For any complex number $z$, the $z$-commutator of two operators
$A, B$ is defined as:
$$[A, B]_z = AB - zBA.$$ The following Lemma holds for $\uS$
 as $\uS'$ and $\uS$ are unitarily equivalent and
  the same properties have been proved for $\uS'$ in [BB].

 \begin{Lemma}
\begin{enumerate}
\item Each monomial $\uS ^{\uk} \omega$ is an eigenvector for $\sum
S_i^*S_i - I$, so that it is a diagonal operator on the standard
basis. In fact,
$$ \sum_{i=1}^n S_i^* S_i (\uS^\uk \omega) =
\left(\sum_{i=1}^n \frac{\|\uS^{\uk + \ue_i}\omega \|^2}
{\|\uS^\uk \omega \|^2} \right) \uS^\uk \omega .$$ Also $\sum
S_i^*S_i - I$ is compact.

\item The commutator $[S_i^* , S_i]$ is as follows:
$$[S_i^* , S_i] \uS^\uk \omega =  \left( \frac{\|\uS^{\uk + \ue_i}\omega\|^2
}{\|\uS^\uk \omega \|^2} - \frac{\|\uS^\uk \omega \|^2}
{\|\uS^{\uk - \ue_i} \omega \|^2} \right) \uS^\uk \omega, \mbox{
when } k_i \neq 0.$$
 If $k_i = 0$, then $[S_i^* , S_i] \uS^\uk \omega =
S_i^*S_i \uS^\uk \omega= \frac{\|\uS^{\uk + \ue_i} \omega
\|^2}{\|\uS^\uk \omega \|^2} \uS^\uk \omega.$

\item $[S_i^*, S_j]_{q_{ij}}$ is compact for all
$1\leq i, j \leq n$.

\end{enumerate}
\end{Lemma}

The map $U^{m,q}:\CS_m \to \Gamma(\mathbb {C}^n)$ given by
$$ U^{m,q}(\sigma) = U^{m,q}_{\sigma}$$
gives the representation of $\CS_m$ on $\Gamma(\mathbb{C}^n).$ It
is easy to see that for all $q=(q_{ij})_{n \times n}, |q_{ij}|=
1,$ the representations are isomorphic or similar by checking the
characters of the representaions. They have same characters. But
for the representations of permuation groups it follows that they
are unitarily equivalent representations. So there exist unitary
$W^q : \Gamma (\mathbb {C}^n) \to \Gamma (\mathbb {C}^n)$ such
that
\begin{equation}
 W^q P_{\Gamma_S (\mathbb{C}^n)} =P_{\Gamma_q (\mathbb{C}^n)}W^q.
\end{equation}
This $W^q$ is not unique as for $k \in \C$ such that $|k|=1,$ the
operator $k W^q$ is also a unitary which satisfy equation (2.4).
We will give one such $W^q$ explicitely.

For  $m \in \N, y_i \in \Lambda$  define $W^{q,m}$ over
$(\C^n)^{\otimes ^m}$ as
$$W^{q,m} ( e_{y_1}\otimes
 \ldots \otimes e_{y_m})
 =   q^{\sigma^{-1}}(x) e_{y_1}\otimes \ldots
\otimes e_{y_m}.$$ where $x=(x_1,\cdots,x_m)$ is the tuple got by
rearranging $(y_1,\cdots,y_m)$ in nondecreasing order and
$\sigma\in \CS_m$ such that $y_i=x_{\sigma(i)}.$ From Proposition
8 its clear that $q^{\sigma^{-1}}(x)$ does not depend upon the
choice of $\sigma.$ And
\begin{eqnarray*}
W^{m,q} P_{\Gamma _S(\mathbb{C}^n)} (e_{y_1}\otimes \ldots \otimes
e_{y_m}) &=& W^{m,q}(\frac{1}{m!} \sum_{\tau \in \CS_m }
e_{y_{\tau^{-1} (1)}}\otimes \ldots
\otimes e_{y_{\tau^{-1}(m)}})\\
&= &  \frac{1}{m!} \sum_{\tau \in \CS_m } q^{(\tau^{-1}
\sigma)^{-1}}(x) e_{y_{\tau^{-1} (1)}}
\otimes \ldots \otimes e_{y_{\tau^{-1} (m)}}\\
&=&  \frac{1}{m!} \sum_{\tau \in \CS_m } q^{\sigma^{-1} \tau}(x)
e_{y_{\tau^{-1}(1)}}\otimes
\ldots \otimes e_{y_{\tau^{-1} (m)}} \\
&= &  \frac{1}{m!} \sum_{\tau \in \CS_m }
q^{\tau}(x_{\sigma(1)},\cdots,x_{\sigma(m)})q^{ \sigma^{-1}}(x)
e_{y_{\tau^{-1}(1)}}\otimes
\ldots \otimes e_{y_{\tau^{-1} (m)}}\\
&= &  P_{\Gamma _q(\mathbb {C}^n)} q^{\sigma^{-1}}(x)
e_{y_1}\otimes \ldots \otimes
e_{y_m}\\
& = &  P_{\Gamma _q(\mathbb {C}^n)}W^{m,q} (e_{y_1}\otimes \ldots
\otimes e_{y_m}).
\end{eqnarray*}
So, $W^{m,q} P_{\Gamma _S(\mathbb {C}^n)}=P_{\Gamma _q(\mathbb
{C}^n)}W^{m,q}$ and $W^q= \oplus^\infty_{m=0} W^{m,q}$ gives the
required unitary which satisfy equation (2.4)(here $W^{0,q}=I$).
Also note that for $\Gamma_q(\C^n)$ and $\Gamma_{q'}(\C^n)$ we
have unitary $W^{q'}(W^{q})^*$ such that
$$W^{q'}(W^{q})^* P_{\Gamma_q(\mathbb {C}^n)}=
P_{\Gamma_{q'}(\mathbb {C}^n)} W^{q'}(W^{q})^*$$

\end{section}

\begin{section}{Dilation of $q$-Commuting Tuples and the Main Theorem}

\setcounter{equation}{0}

\begin{Definition}
{\em Let $\uT =(T_1, \ldots , T_n)$ be a contractive tuple on a
Hilbert space $\CH .$ The operator $\Delta _{\uT}=
[I-(T_1T_1^*+\cdots +T_nT_n^*)]^{\frac{1}{2}}$ is called the {\em
defect operator \/} of $\uT$ and the subspace $\overline {\Delta
_\uT(\CH)}$ is called the {\em defect space\/} of $\uT .$ The
tuple $\uT $ is said to be {\em pure \/ } if $\sum _{\alpha \in
\Lambda ^m}\uT^{\alpha }(\uT^{\alpha})^*$ converges to zero in
strong operator topology as $m$ tends to infinity. }
\end{Definition}

When $\sum T_iT_i^*=I$, we have $\sum _{\alpha \in \Lambda
^m}\uT^{\alpha }(\uT^{\alpha})^* =I$ for all $m$ and hence $\uT$
is not pure.  Let $\uT $ be a pure contractive tuple on $\CH .$
Take $\tilde {\CH}=\Gamma (\mathbb{C}^n)\otimes \overline{\Delta
_{\uT}(\CH)},$ and define an operator $A:\CH \to \tilde {\CH}$ by
\begin{equation}
Ah= \sum _{\alpha }e^{\alpha }\otimes \Delta _{\uT}(\uT ^{\alpha
})^*h,
\end{equation}
where the sum is taken over all $\alpha \in \tilde {\Lambda}$
(this operator was used by  Popescu and Arveson in [Po3], [Po4],
[Ar2] and for $q$-commuting case by Bhat and Bhattacharyya in
[BB]). $A$ is an isometry  and we have $\uT ^\alpha =A^*(\uV
^\alpha \otimes I)A$ for all $\alpha \in {\tilde \Lambda}$ (see
[Po4]). Also the tuple ${\tilde \uV}= (V_1\otimes I, \ldots ,
V_n\otimes I)$ of operators on ${\tilde \CH }$ is a realization of
the minimal noncommuting dilation of $\uT$.

Let $C^*(\uV)$, and  $C^*(\uS)$ be unital $C^*$-algebras generated
by tuples $\uV$ and $\uS$ (defined in the Introduction) on Fock spaces $\Gamma
({\mathbb{C}}^n)$ and $ \Gamma _q({\mathbb{C}}^n)$ respectively.
 For any $\alpha , \beta \in
{\tilde \Lambda }$, $\uV ^{\alpha }(I-\sum V_iV_i^*) (\uV ^{\beta
})^*$ is the rank one operator $x\mapsto  \langle e^{\beta },
x\rangle e^\alpha ,$ formed by basis vectors $e^{\alpha },
e^{\beta }$ and so $C^*(\uV)$ contains all compact operators. Similarly
  we see that $C^*(\uS)$ also contains all compact
operators of $\Gamma _q({\mathbb{C}}^n)$. As $V_i^*V_j=\delta
_{ij}I,$ it is easy to see that $C^*(\uV)= ~~\overline {\mbox
{span}}~~\{\uV^\alpha (\uV^\beta)^*: \alpha , \beta \in {\tilde
\Lambda}\}.$  As $q_{ij}$-commutators $[S_i^*, S_j]_{q_{ij}}$ are
compact for all $i,j$, we can also get
$C^*(\uS)=~~\overline {\mbox {span}}~~\{\uS^\alpha (\uS^\beta)^*:
\alpha , \beta \in {\tilde \Lambda}\}.$

 Consider a contractive tuple $\uT $ on a Hilbert space $\CH$.
For $0< r <1 $ the tuple $r\uT= (rT_1, \ldots , rT_n)$ is clearly
a pure contraction. So by equation (2.4) we have an isometry $A_r: \CH \to
\Gamma ({\mathbb{C}}^n)\otimes \overline {\Delta _r(\CH)}$ defined
by
$$A_rh = \sum _{\alpha }e^{\alpha }\otimes \Delta _r((r\uT) ^\alpha )^*h, ~~h\in \CH ,$$
where $\Delta _r=(I-r^2\sum T_iT_i^*)^{\frac{1}{2}}.$ So for every
$0<r<1$ we have a completely positive map $\psi _r: C^*(\uV)\to
\CB(\CH)$ defined by
$\psi _r(X)=A_r^*(X\otimes I)A_r, ~~X\in C^*(\uV).$
By taking limit as $r$ increases to 1 (See [Po1-4] for details),
we get a unital completely positive map $\psi $ from $C^*(\uV)$
to $\CB(\CH)$ (Popescu's Poisson transform) satisfying
$$\psi (\uV^\alpha (\uV^\beta)^*)= \uT^\alpha (\uT^\beta )^* ~~\mbox{for} ~ \alpha , \beta \in
{\tilde \Lambda}.$$ As $C^*(\uV)= ~~\overline {\mbox
{span}}~~\{\uV^\alpha (\uV^\beta)^*: \alpha , \beta \in {\tilde
\Lambda}\},$ $\psi $ is the unique such completely positive map.
Let the minimal Stinespring dilation of $\psi $ be unital
$*$-homomorphism $\pi : C^*(\uV)\to \CB({\tilde \CH })$ where
 ${\tilde \CH}$ is a  a Hilbert space containing $\CH $, and
$$\psi (X)=P_{\CH}\pi(X)|_{\CH} ~~\forall X\in C^*(\uV),$$
and $\overline {\mbox {span}}~\{ \pi(X)h: X\in C^*(\uV), h\in
\CH\}= {\tilde \CH}.$ Let ${\tilde \uV}=({\tilde V_1}, \ldots ,
{\tilde V_n})$ where ${\tilde V_i}=  \pi(V_i)$ and so ${\tilde \uV}$
is the unique standard noncommuting dilation of $\uT$ and clearly
 $\tilde {(V_i)}^*$ leaves $\CH $ invariant. If $\uT$ is $q$-commuting,
by considering $C^*(\uS)$ instead of
$C^*(\uV)$, and restricting $A_r$ in the range to $\Gamma
_q({\mathbb{C}}^n)$, and taking limits as $r$ increases to 1 as before  we
would get the unique unital completely positive map $\phi :
C^*(\uS)\to \CB (\CH),$ (also see [BB]) satisfying

\begin{equation}
\phi (\uS^\alpha (\uS^\beta)^*)= \uT^\alpha (\uT^\beta )^* ~~~~~
\alpha , \beta \in {\tilde \Lambda}.
\end{equation}

\begin{Definition}
{\em  Let $\uT$ be a $q$-commuting tuple. Then we have a unique
unital completely positive map $\phi : C^*(\uS)\to \CB (\CH)$
satisfying equation (3.2). Consider the minimal Stinespring
dilation of $\phi .$ Here we have a Hilbert space $\CH _1 $
containing $\CH$ and a unital $*$-homomorphism $\pi _1 :
C^*(\uS)\to \CB(\CH _1),$ such that
$$ \phi (X)=P_{\CH}\pi _1(X)|_{\CH} ~~~~\forall X\in C^*(\uS),$$
and $\overline {\mbox {span}}~\{ \pi _1(X)h: X\in C^*(\uS), h\in
\CH \}= { \CH _1}.$ Let $\tilde{S_i}=\pi_1(S_i)$ and
$\tilde{\uS}=(\tilde{S_1},\ldots,\tilde{S_n}).$ Then $\tilde
{\uS}$ is called the {\em standard $q$-commuting dilation \/} of
$\uT. $}
\end{Definition}

 Standard $q$-commuting dilation is also unique up to
unitary equivalence as minimal Stinespring dilation is unique up to unitary
equivalence.

\begin{Lemma} Suppose $\uT=(T_1, \cdots , T_n)$ is a $q$-commuting tuple on a
Hilbert Space $\CH$ and let $A$ be the operator introduced in
Equation (3.1). Then there exist a Hilbert space $K$ such that
$(S_1\otimes I_{\CK}, \ldots , S_n\otimes I_{\CK})$ is a dilation
of $\uT$ and dim $(\CK)=$ rank $~(\Delta _{\uT}).$
\end{Lemma}

\noindent {\sc Proof:} $A(h)=\sum_\alpha e^\alpha \otimes \Delta_\uT (\uT^\alpha)^*h $ for $h \in \CH$
where the sum is over $\alpha \in \tilde{\Lambda}.$
For a given $\uk=(k_1,\cdots,k_n)$ such that $|\uk|=m$ let us denote
$e_1^{k_1} \otimes \cdots \otimes e_n^{k_n}$ by $e_{x_1}\otimes \cdots \otimes
e_{x_m}, 1 \leq x_m \leq n$ in the following calculation.
\begin{eqnarray*}
A(h)&=& \sum_{m=0}^\infty \sum_{\sigma \in \CS_m}
e_{x_{\sigma^{-1}(1)}} \cdots e_{x_{\sigma^{-1}(m)}} \otimes
\Delta_\uT ( T_{x_{\sigma^{-1}(1)}} \cdots T_{x_{\sigma^{-1}(m)}})^*h\\
&=& \sum_{m=0}^\infty \sum_{\sigma \in \CS_m}
e_{x_{\sigma^{-1}(1)}} \cdots e_{x_{\sigma^{-1}(m)}} \otimes
\Delta_\uT  \overline{(q^\sigma (x))^{-1}}( T_{x_1} \cdots T_{x_m})^*h \\
&=& \sum_{m=0}^\infty \sum_{\sigma \in \CS_m} q^\sigma (x)
e_{x_{\sigma^{-1}(1)}} \cdots e_{x_{\sigma^{-1}(m)}} \otimes
\Delta_\uT   ( T_{x_1} \cdots T_{x_m})^*h \\
&=& \sum_{m=0}^\infty  (m!) P_m e_{x_1} \cdots e_{x_m} \otimes
\Delta_\uT   ( T_{x_1} \cdots T_{x_m})^*h
\end{eqnarray*}

So the range of $A$ is contained in ${\tilde
\CH}_q =  \Gamma _q({\mathbb{C}}^n)\otimes \overline {\Delta
_{\uT}(\CH)}$. In other words now $\CH $ can be considered as a
subspace of ${\tilde \CH}_q$. Moreover,   $\tilde
{\uS}=(S_1\otimes I, \ldots , S_n\otimes I)$, as a tuple of
operators in ${\tilde \CH}_q$ is the standard $q$-commuting
dilation of $(T_1, \ldots T_n).$ More abstractly we can get a
Hilbert space $\CK$ such that  $\CH$ can be isometrically embedded
in $\Gamma _q(\mathbb{C}^n)\otimes \CK$ and $(S_1\otimes I_{\CK},
\ldots , S_n\otimes I_{\CK})$ is a dilation of $\uT$ and
$\overline {\mbox {span}}\{(\uS ^\alpha\otimes I_{\CK})h: h\in
\CH, \alpha \in \tilde {\Lambda }\}=\Gamma _q(\mathbb{C}^n)\otimes
\CK$.  There is a unique such dilation and up to unitary equivalence
and dim $(\CK)=$ rank $~(\Delta _{\uT}).$ \qed

\begin{Theorem}

Let $\uT $ be a pure contractive tuple on a Hilbert space $\CH$.
\begin{enumerate}
\item Then the maximal $q$-commuting piece ${\tilde \uV}^q$ of the
standard noncommuting dilation ${\tilde \uV}$ of $\uT$ is a
realization of the
 standard $q$-commuting dilation of $\uT ^q$ if and only
if $\overline {\Delta _\uT(\CH)}=\overline {\Delta _\uT(\CH
^q(\uT))}.$ And if $\overline {\Delta _\uT(\CH)}=\overline {\Delta _\uT(\CH
^q(\uT))}$ then rank $(\Delta _{\uT})=$ rank $(\Delta
_{\uT^q})=$ rank $(\Delta _{{\tilde \uV}})=$ rank $(\Delta
_{{\tilde \uV}^q}).$

\item Let the  standard noncommuting dilation of $\uT $ be ${\tilde \uV}$.
If rank $\Delta _{\uT}$
 and rank $\Delta _{\uT ^q}$ are finite and equal then ${\tilde \uV}^q$
 is a realization of the  standard $q$-commuting  dilation of $\uT ^q$.
 \end{enumerate}
\end{Theorem}

\noindent {\sc Proof:}
The proof is similar to the proofs  of that of Theorem 10 and Remark 11
of [BBD].  \qed

If  the ranks of both $\Delta _{\uT}$
 and  $\Delta _{\uT ^q}$ are infinite then we can not ensure that
$\overline {\Delta _\uT(\CH)}=\overline {\Delta _\uT(\CH
^q(\uT))}$ and hence can not ensure the converse of the last Theorem,
as seen by the following example.
For any $n \geq 2$ consider the Hilbert space $\CH _0=\Gamma_q
(\mathbb {C}^n) \otimes \CM$ where $\CM$ is of infinite dimension
and let $\uR=(S_1 \otimes I,\cdots, S_n \otimes I)$ be a $q$-commuting
pure contractive $n$-tuple. Infact one can take any $\uR$ to be any  $q$-commuting
pure $n$-tuple on some Hilbert space $\CH_0$ with
 $\overline {\Delta _\uR(\CH _0)}$ of  infinite dimensional. Suppose
 $P_k=(p^k_{ij})_{n \times n}$ for $1 \leq k \leq n$ are $n \times n$ matrices
with complex entries such that
$$p^k_{ij}= \left \{
\begin{array}{ccc}
t_k  & \mbox {~if~} i=k,j=k+1 \\
0 & \mbox{~otherwise~}
\end{array} \right . \mbox{~for~} 1 \leq k < n \mbox{~and~}
p^n_{ij}= \left\{ \begin{array}{ccc}
t_n  & \mbox {~if~}i=n,j=1\\
0 & \mbox{~otherwise~}
\end{array} \right .$$
where $t_k$'s are complex numbers satisfying $0< |t_k| <1.$
Let $\CH =\CH _0\oplus \mathbb{C}^n$. Take $\uT=(T_1, \cdots, T_n)$ where $T_k$
for $1 \leq k \leq n$  be operators on $\CH $ defined by
$$T_k=\left[\begin{array}{ccc}
R_k & \\
& P_k
\end{array}\right] \mbox{~for~} 1\leq k \leq n.$$
 So $\uT$
is a pure contractive tuple, the maximal $q$-commuting piece of
$\uT$ is $\uR$ and  $\CH ^q(\uT)=\CH _0$ (by Corollary 7). Here
rank$~(\Delta _{\uT^q})=$ rank$~(\Delta _{\uT})=\infty$ but
$\overline {\Delta _\uT(\CH)}=\overline {\Delta _\uR(\CH
_0)}\oplus \mathbb{C}^n.$ But the converse of Theorem 18 holds
when rank of $\Delta_{\uT}$ is finite.

 Consider the case when $\uT $ is a
$q$-commuting tuple on Hilbert space $\CH $ satisfying $\sum
T_iT_i^* =I$.
 As $C^*(\uS)$ contains the ideal of all
compact operators by standard $C^*$-algebra theory we have a
direct sum decomposition of $\pi _1$ as follows. Take
${\mathcal H}_1 = {\mathcal H}_{1C} \oplus {\mathcal H}_{1N}$
where ${\mathcal H}_{1C} = \overline{\mbox{span}}\{\pi_1(X)h :
h\in {\mathcal H}, X \in C^*(\underline{S})$ and $X$ is
compact$\}$ and ${\mathcal H}_{1N}$ is the orthogonal complement of it.
Clearly ${\mathcal H}_{1C}$ is a reducing subspace for $\pi _1$.
Therefore  $\pi _1=\pi _{1C}\oplus \pi _{1N}$ where $\pi_{1C}(X) =
P_{{\mathcal H}_{1C}} \pi _1(X) P_{{\mathcal H}_{1C}}$,
$\pi_{1N}(X) = P_{{\mathcal H}_{1N}} \pi _1(X) P_{{\mathcal
H}_{1N}}$. Also $\pi _{1C}(X)$ is just the identity
representation with some multiplicity. Infact  $\CH _{1C}$
can be written as $\CH _{1C}= \Gamma _q({\mathbb{C}}^n)\otimes
\overline {\Delta _{\uT}(\CH)}$ (see Theorem 4.5 of [BB]) and $\pi
_{1N}(X)=0$ for compact $X$. But $\Delta _{\uT}(\CH)=0$ and commutators
$[S_i^*,S_i]$ are compact. So if we take $W_i=\pi_{1N}(S_i)$,
$\underline {W}=(W_1, \ldots , W_n)$ is a tuple of normal operators.
 It follows that the standard $q$-commuting dilation of $\uT $ is a tuple of normal operators.

\begin{Definition}
{\em A $q$-commuting $n$-tuple $\uT =(T_1, \ldots , T_n) $ of operators on a
Hilbert space $\CH $ is called a $q$-{spherical unitary \/} if
 each $T_i$ is normal and $T_1T_1^*+\cdots
+T_nT_n^*=I.$}
\end{Definition}

If $\CH $ is a finite dimensional Hilbert space and $\uT $ is a
$q$-commuting tuple on $\CH $ satisfying $\sum T_iT_i^* =I$, then
 $\uT $ a spherical unitary because each $T_i$ would be subnormal
 and all finite dimensional subnormal operators are normal (see [Ha]).

\begin{Theorem} (Main Theorem)
Let $\uT$ is a $q$-commuting contractive tuple on a Hilbert
space $\CH .$ Then the maximal $q$-commuting piece of the standard
noncommuting dilation of $\uT $ is a realization of the standard
$q$-commuting dilation of $\uT$.
\end{Theorem}

\noindent {\sc Proof of the theorem 19:} Let $\tilde{\uS}$ denote
the standard $q$-commuting dilation of $\uT $ on a Hilbert space
$\CH _1$ and we follow the notations as in section 2. As  $\uS $
is also a contractive tuple, we have a unique unital completely
positive map $\eta : C^*(\uV)\to C^*(\uS)$, satisfying
$$\eta (\uV^\alpha (\uV^\beta)^*)= \uS^\alpha (\uS^\beta )^* ~~~~ \alpha , \beta \in
{\tilde \Lambda}.$$ It is easy to see that $\psi = \phi \circ \eta
$. Let unital $*$-homomorphism $\pi _2 : C^*(\uV)\to \CB(\CH _2)$
for some  Hilbert space $\CH _2 $ containing  $\CH _1,$ be the
minimal Stinespring dilation of the map $\pi _1\circ \eta :
C^*(\uV)\to \CB(\CH _1)$   such that $\pi _1\circ \eta (X)=P_{\CH
_1 }\pi _2(X)|_{\CH _1}, ~~ ~~\forall X\in C^*(\uV),$ and
$\overline {\mbox {span}}~\{ \pi _2(X)h: X\in C^*(\uV), h\in \CH
_1\}= {\CH _2}.$ We get the following commuting diagram.

\hskip1.5in \begin{picture}(200,115)

\put(0,20){$C^*(\uV)$}

\put(40,21){$\longrightarrow$}

\put(70,20){$C^*(\uS)$}

\put(110,21){$\longrightarrow$}

\put(140,20){$\mathcal{B}(\CH)$}

\put(140,60){$\mathcal{B}(\CH_1)$}

\put(140,100){$\mathcal{B} (\CH_2)$}

\put(40,30){\vector(4,3){80}}

\put(110,30){\vector(4,3){20}}

\put(150,40){$\downarrow $}

\put(150,80){$\downarrow $}

\put(45,10){$\eta$}

\put(115,10){$\phi$}

\put(110,40){$\pi_1$}

\put(70,70){$\pi_2$}

\end{picture}

\noindent where all the down arrows are compression maps,
horizontal arrows are unital completely positive maps and diagonal
arrows are unital $*$-homomorphisms. Let ${\hat  \uV}=({\hat V_1},
\ldots , {\hat V_n})$ where $ {\hat V_i}=\pi _2(V_i).$ We would
show that ${\hat \uV}$ is the standard noncommuting dilation of
$\uT.$  We have this result if we can show that $\pi _2$ is a
minimal dilation of $\psi=\phi \circ\eta $  as  minimal
Stinespring dilation is unique up to unitary equivalence. For this
first we show that  ${\tilde {\uS }} =(\pi _1(S_1), \ldots , \pi
_1(S_n))$ is the maximal $q$-commuting piece of $\hat {\uV}$.

 First we consider a
particular case when $\uT $ is a $q$-spherical unitary on a Hilbert space $\CH .$
In this case we would show that  standard commuting dilation and
the maximal $q$-commuting piece of the standard noncommuting
dilation of $\uT $ is itself.

We have $\phi(\uS^\alpha (I-\sum S_iS_i^*)(\uS ^\beta )^*) =\uT ^\alpha
(I-\sum T_iT_i^*)(\uT ^\beta)^*=0$ for any $\alpha , \beta \in
{\tilde \Lambda}.$ This forces that $\phi (X)=0$ for any compact
operator $X$ in $C^*(\uS).$ Now as the $q_{ij}$-commutators
$[S_i^*, S_j]_{q_{ij}}$ are all compact we see that $\phi $ is a
unital $*$-homomorphism. So the minimal Stinespring dilation of
$\phi $ is itself and standard commuting dilation of $\uT $ is itself.
Next we would show that
the maximal $q$-commuting piece of the standard noncommuting
dilation of $\uT $ is itself. The presentation of the standard
noncommuting dilation which we would use is taken from [Po1]. The dilation space
$\tilde {\CH} $ can be decomposed as $ \tilde {\CH } ={\mathcal H} \oplus
(\Gamma({\mathbb{C}}^n)\otimes\mathcal{D})$ where $\mathcal{D}$ is
the closure of the range of  operator
\[ D:\underbrace{{\mathcal H}\oplus \cdots \oplus {\mathcal H}}_{n~
copies} \rightarrow \underbrace{{\mathcal H} \oplus \cdots \oplus
{\mathcal H}}_{n ~copies} \] and $D$ is the positive square root
of
\[ D^2=[\delta_{ij}I-T_i^*T_j]_{n \times n}. \]
For convenience, at some places we would identify
$\underbrace{{\mathcal H}\oplus \cdots \oplus {\mathcal
H}}_{n~copies}$ with ${\mathbb{C}}^n\otimes {\mathcal H}$ so that
 $(h_1, \ldots ,h_n)= \sum_{i=1}^n e_i \otimes h_i.$
Then
\begin{equation}
 D(h_1, \ldots ,h_n)= D(\sum_{i=1}^n e_i \otimes h_i)=\sum_{i=1}^n
e_i\otimes (h_i - \sum _{j=1}^nT^*_iT_j h_j)
\end{equation}
and the standard noncommuting dilation $\tilde{V_i}$
\begin{equation}
 \tilde{V_i}(h\oplus \sum_{\alpha \in
{\tilde{\Lambda}}}e^{\alpha}\otimes d_{\alpha})= T_ih\oplus
D(e_i\otimes h )\oplus e_i\otimes(\sum_{\alpha \in
{\tilde{\Lambda}}} e^{\alpha}\otimes d_{\alpha})
\end{equation}
for $h \in \CH$,
 $d_{\alpha} \in
\mathcal{D}$ for $\alpha \in \tilde{\Lambda}$, and $1 \leq i \leq
n$ ($\mathbb{C}^n\omega \otimes \CD$ has been identified with
$\CD$).

 We have $$T_iT_i^*=
T_i^*T_i \mbox {~and~} T_jT_i= q_{ij}T_iT_j \forall 1 \leq i, j
\leq n .$$ Also by Fuglede-Putnam Theorem ([Ha] [Pu])
$$T_j^*T_i=\oq_{ij} T_iT_j^* = q_{ji} T_i T_j^* \mbox {~and~}
 T_j^*T_i^*= q_{ij}T_i^*T_j^* \forall 1 \leq i, j \leq n .$$
 As $\sum T_iT_i^*=I$, by direct computation
 $D^2$ is seen to be a projection. So, $D=D^2$. Note that $q_{ij}\oq_{ij} = 1,$
i.\@ e., $\oq_{ij} = q_{ji}.$ Then
we get
\begin{eqnarray}
 D(h_1, \ldots ,h_n)
&= &  \sum_{i,j=1}^n e_i \otimes T_j (T_j^*h_i - \oq_{ji}T_i^*h_j)
 = \sum_{i,j=1}^n e_i \otimes T_j (h_{ij})
\end{eqnarray}
where $h_{ij} = T_j^*h_i - \oq_{ji} T_i^*h_j = T_j^*h_i - q_{ij}
T_i^*h_j$ for $1 \leq i, j \leq n$. Note that $h_{ii}=0$ and
$h_{ji}=-\oq_{ij}h_{ij}.$

As clearly $\CH \subseteq \tilde{\CH}^q(\uV),$ lets begin with
$y\in \CH ^{\perp}\bigcap \tilde {\CH }^q(\tilde {\uV}). $ We wish
 to show that $y=0$. Decompose $y$ as $y=0\oplus \sum _{\alpha \in \tilde {\Lambda }}
 e^{\alpha }\otimes y_{\alpha },$ with $y_{\alpha }\in \CD .$ We assume $y\neq 0$
and arrive at a contradiction. If for some $\alpha $, $y_{\alpha }\neq 0$, then
 $\langle \omega \otimes y_{\alpha }, (\tilde {\uV}^{\alpha })^*y\rangle
 =\langle e^{\alpha }\otimes y_{\alpha }, y\rangle =
 \langle y_{\alpha }, y_{\alpha }\rangle \neq  0.$
 Since  $(\tilde {\uV}^{\alpha })^*y\in \tilde {\CH }^q(\tilde {\uV}),$
  we can assume $\|y_0\|=1$ without loss of generality.
Taking $\tilde {y}_m =\sum _{\alpha \in \Lambda ^m}e^{\alpha
}\otimes y_{\alpha}$, we get $y=0\oplus \oplus _{m\geq 0}\tilde
{y}_m.$ $D$ being a projection its range is closed and as $y_0 \in
{\mathcal D}$, there exist some $(h_1, \ldots ,h_n)$ such that
$y_0=D(h_1, \ldots ,h_n)$. Let $\tilde{x_0}=\tilde{y_0}=y_0$,
$\tilde{x_1}= \sum_{i,j=1}^n e_i \otimes D(e_j\otimes h_{ij}),$
and for $m\geq 1$,
$$ \tilde{x}_m = \sum_{i_1, \ldots ,i_{m-1},i,j=1}^n e_{i_1}
\otimes \cdots \otimes e_{i_{m-1}} \otimes e_i \otimes D(e_j
\otimes (\prod_{1 \leq r < s \leq
m-1}q_{i_ri_s})(\prod_{k=1}^{m-1}q_{i_ki}q_{i_kj}) T_{i_1}^*
\ldots T_{i_{m-1}}^* h_{ij}).$$
 So $\tilde{x}_m \in ({\mathbb
C}^n)^{\otimes m} \otimes {\mathcal D}$ for all $m \in
{\mathbb{N}}$. As $\uT$ is  $q$-commutating $n$-tuple and $D$ is a
projection, we have

\begin{eqnarray*}
\sum_{1\leq i < j \leq n} (q_{ij}\tilde{V}_i \tilde{V}_j -
\tilde{V}_j \tilde{V}_i)q_{ji}h_{ij}
&=&\sum _{1\leq i<j\leq n}(q_{ij}T_iT_j-T_jT_i)q_{ji}h_{ij}\\
 & &+\sum_{1\leq i < j \leq n}D(e_i\otimes T_j h_{ij} -
q_{ji}e_j\otimes T_i h_{ij})\\
& &
  +\sum_{1\leq i < j \leq n}(e_i\otimes D(e_j\otimes h_{ij})
   -q_{ji}e_j\otimes D(e_i \otimes h_{ij}))\\
& =& 0 + D(\sum _{i,j = 1}^n e_i\otimes T_jh_{ij}) + \sum_{ i,j=1}^n e_i\otimes D(e_j\otimes h_{ij})\\
& =& D^2(h_1, \ldots ,h_n) + \sum_{ i,j=1}^n e_i\otimes
D(e_j\otimes h_{ij})\\
&=& \tilde{x}_0 +\tilde{x}_1.
\end{eqnarray*}
So by Proposition 6, $\langle y, \tilde{x}_0 +\tilde{x}_1 \rangle
=0$ . Next let $m \geq 2.$
\begin{eqnarray*}
& & \sum_{i_1, \ldots ,i_{m-1}=1}^n \tilde{V}_{i_1} \ldots
\tilde{V}_{i_{m-1}}\{\sum_{i,j=1}^n (q_{ij}
\tilde{V}_i\tilde{V}_j-\tilde{V}_j\tilde{V}_i)(\prod_{1 \leq r < s
\leq m-1}q_{i_ri_s})(\prod_{k=1}^{m-2}q_{i_kj})(T_i^* T_{i_1}^*
\ldots T_{i_{m-2}}^* h_{i_{m-1}j})\}\\
& =&  \sum_{i_1, \ldots ,i_{m-1}=1}^n e_{i_1} \otimes \ldots
\otimes e_{i_{m-1}}\otimes [\sum_{i,j=1}^n D((\prod_{1 \leq r < s
\leq m-1} q_{i_ri_s})(\prod_{k=1}^{m-2}q_{i_kj})(q_{ij} e_i\otimes
T_jT_i^*T_{i_1}^*\ldots T_{i_{m-2}}^*h_{i_{m-1}j}\\
& & - e_j\otimes T_iT_i^*T_{i_1}^*\ldots
T_{i_{m-2}}^*h_{i_{m-1}j})) +   \sum_{i,j=1}^n(\prod_{1 \leq r < s
\leq m-1}\ q_{i_ri_s})(\prod_{k=1}^{m-2}q_{i_kj}) \{q_{ij} e_i
\otimes D(e_j \otimes \\
& & T_i^*T_{i_1} ^* \ldots T_{i_{m-2}}^*h_{i_{m-1}j}) - e_j\otimes
D(e_i \otimes T_i^*T_{i_1}^*\ldots
T_{i_{m-2}}^*h_{i_{m-1}j})\}]\\
 &= &- \sum_{i_1,\ldots
,i_{m-1}=1}^n e_{i_1} \otimes \cdots \otimes e_{i_{m-1}} \otimes
\{(\sum_{j=1}^n (\prod_{1 \leq r < s \leq m-1}
q_{i_ri_s})(\prod_{k=1} ^{m-2}q_{i_kj}) D(e_j \otimes
T_{i_1}^* \ldots T_{i_{m-2}}^* h_{i_{m-1}j})\}\\
& & + \sum_{i_1,\ldots ,i_{m-1}=1}^n e_{i_1} \otimes \ldots
\otimes e_{i_{m-1}} \otimes \{\sum_{i,j=1}^n e_i \otimes D(e_j
\otimes q_{ij} (\prod_{1 \leq r
< s \leq m-1}q_{i_ri_s})(\prod_{k=1}^{m-2} q_{i_kj})\\
& & (T_i^*T_{i_1}^* \ldots T_{i_{m-2}}^* h_{i_{m-1}j}))
-\sum_{i,j=1}^n e_i \otimes D(e_j\otimes (\prod_{1 \leq r < s \leq
m-1}q_{i_ri_s})(\prod_{k=1}^{m-2} q_{i_ki})
 (T_j^*T_{i_1}^* \ldots T_{i_{m-2}}^* h_{i_{m-1}i}))\}\\
& & \mbox{(in the term above, $i$ and $j$ have been interchanged
in the last summation)}\\
 &= &- \sum_{i_1,\ldots ,i_{m-1}=1}^n
e_{i_1} \otimes \cdots \otimes e_{i_{m-1}} \otimes \{\sum_{j=1}^n
(\prod_{1 \leq r < s \leq m-2} q_{i_ri_s})(\prod_{k=1}
^{m-2}q_{i_ki}q_{i_kj}) D(e_j \otimes
T_{i_1}^* \ldots T_{i_{m-2}}^* h_{ij})\}\\
& & + \sum_{i_1,\ldots ,i_{m-1}=1}^n e_{i_1} \otimes \ldots
\otimes e_{i_{m-1}}\otimes \{\sum_{i,j=1}^n e_i \otimes D(e_j\\
& & \otimes  (\prod_{1 \leq r < s \leq m-1} q_{i_ri_s})
q_{ij}(\prod_{k=1}^{m-2} q_{i_kj}) (T_i^*T_{i_1}^* \ldots
T_{i_{m-2}}^*T_j^* h_{i_{m-1}}
-q_{i_{m-1}j}T_i^*T_{i_1}^* \ldots T_{i_{m-2}}^*T_{i_{m-1}}^* h_j )\\
& & -(\prod_{1 \leq r < s \leq m-2}q_{i_ri_s})(\prod_{k=1}^{m-2}
q_{i_ki}) (T_j^*T_{i_1}^* \ldots T_{i_{m-2}}^* T_i^*h_{i_{m-1}}-
q_{i_{m-1}i}T_j^*T_{i_1}^* \ldots T_{i_{m-2}}^* T_{i_{m-1}}^*h_i ))\}\\
& & \mbox{(in the  term above, index $i_{m-1}$ has been replaced
by $i$ in the first summation)}\\
& =& - \sum_{i_1,\ldots ,i_{m-2},i,j=1}^n e_{i_1} \otimes \cdots
\otimes e_{i_{m-2}} \otimes e_i \otimes  (\prod_{1 \leq r < s \leq
m-2} q_{i_ri_s})(\prod_{k=1}^{m-2}q_{i_ki}q_{i_kj})
D(e_j \otimes T_{i_1}^* \ldots   T_{i_{m-2}}^* h_{ij})\\
& & + \sum_{i_1,\ldots  ,i_{m-1},i,j=1}^n e_{i_1} \otimes \cdots
\otimes e_{i_{m-1}} \otimes e_i \otimes (\prod_{1 \leq r < s \leq
m-1}q_{i_ri_s})(\prod_{k=1}^{m-1}q_{i_ki}q_{i_kj})D(e_j\otimes
T_{i_1}^* \ldots   T_{i_{m-1}}^* h_{ij})\\
 & = & -\tilde{x}_{m-1} +
\tilde{x}_m.
\end{eqnarray*}
Hence by proposition 6, $\langle y, \tilde{x}_{m-1} - \tilde{x}_m\rangle =0.$
Further we compute $\|\tilde{x}_{m}\| $ for all $m \in
{\mathbb{N}}$.
\begin{eqnarray*}
\|\tilde{x}_{m}\|^2 & =& \langle \sum_{i_1,\ldots
,i_{m-1},i,j=1}^n e_{i_1} \otimes \cdots  \otimes e_{i_{m-1}}
\otimes e_i \otimes D(e_j\otimes (\prod_{1 \leq r < s \leq
m-1}q_{i_ri_s})(\prod_{k=1}^{m-1}q_{i_ki}q_{i_kj})T_{i_1}^* \ldots
T_{i_{m-1}}^* h_{ij}),\\
& & \sum_{i_1',\ldots ,i_{m-1}',i',j'=1}^n e_{i_1'} \otimes \cdots
\otimes e_{i_{m-1}'} \otimes e_{i'} \otimes D(e_{j'}\otimes
(\prod_{1 \leq r' < s' \leq
m-1}q_{i_{r'}i_{s'}})\\
& & (\prod_{k'=1}^{m-1}q_{i'_{k'}i'}
 q_{i'_{k'}j'})T_{i_1'}^*
\ldots   T_{i_{m-1}'}^* h_{i'j'})\rangle\\
&= &  \sum_{i_1,\ldots  ,i_{m-1},i=1}^n \langle \sum_{j=1}^n
D(e_j\otimes (\prod_{1 \leq r < s \leq m-1}q_{i_ri_s}
)(\prod_{k=1}^{m-1}q_{i_ki}q_{i_kj})T_{i_1}^* \ldots
T_{i_{m-1}}^* h_{ij}),\\
& & \sum_{j'=1}^n D(e_{j'}\otimes (\prod_{1 \leq r' < s' \leq
m-1}q_{i_{r'}i_{s'}})(\prod_{k'=1}^{m-1}
q_{i_{k'}i}q_{i_{k'}j'})T_{i_1}^* \ldots T_{i_{m-1}}^*
h_{ij'})\rangle\\ &= &  \sum_{i_1,\ldots ,i_{m-1},i=1}^n \langle
D(\sum_{j=1}^n e_j\otimes (\prod_{1 \leq r < s \leq m-1}q_{i_ri_s}
)(\prod_{k=1}^{m-1}q_{i_ki}q_{i_kj})T_{i_1}^* \ldots T_{i_{m-1}}^*
h_{ij}),\\
& &  \sum_{j'=1}^n e_{j'}\otimes (\prod_{1 \leq r' < s'\leq
m-1}q_{i_{r'}i_{s'}}
)(\prod_{k'=1}^{m-1}q_{i_{k'}i}q_{i_{k'}j'})T_{i_1}^*
\ldots T_{i_{m-1}}^*  h_{ij'}\rangle\\
 & =&
\sum_{i_1,..,i_{m-1},i=1}^n \langle ( \prod_{1 \leq r < s \leq
m-1}q_{i_ri_s} )(\prod_{k=1}^{m-1}q_{i_ki}q_{i_kj})\{
\sum_{j,l=1}^n e_j\otimes T_l (T_l^*T_{i_1}^* \ldots T_{i_{m-1}}^*
h_{ij}- \\
& & q_{jl}T_j^* T_{i_1}^*\ldots T_{i_{m-1}}^* h_{il})\},
\sum_{j'=1}^n (\prod_{1 \leq l' < s'\leq m-1}q_{i_{l'}i_{s'}}
)(\prod_{k'=1}^{m-1}q_{i_{k'}i}q_{i_{k'}j'}) e_{j'}\otimes
T_{i_1}^* \ldots   T_{i_{m-1}}^* h_{ij'}\rangle\\
& =& \sum_{i_1,...,i_{m-1},i,j=1}^n \langle ( \prod_{1 \leq r < s
\leq m-1}q_{i_ri_s}
)(\prod_{k=1}^{m-1}q_{i_ki}q_{i_kj})\sum_{l=1}^n T_l
(T_l^*T_{i_1}^* \ldots   T_{i_{m-1}}^* h_{ij}- q_{jl}T_j^*
T_{i_1}^*\ldots T_{i_{m-1}}^* h_{il}), \\
& & (\prod_{1 \leq r' < s'\leq m-1}q_{i_{r'}i_{s'}}
)(\prod_{k'=1}^{m-1}q_{i_{k'}i}q_{i_{k'}j})
T_{i_1}^* \ldots   T_{i_{m-1}}^* h_{ij}\rangle\\
& =& \sum_{i,j=1}^n \langle h_{ij},h_{ij}\rangle -
\sum_{i_1,...,i_{m-1},i,j,l=1}^n \langle T_{i_{m-1}} \ldots
T_{i_1} T_j^* T_l T_{i_1}^* \ldots   T_{i_{m-1}}^*
h_{il}),h_{ij}\rangle
\end{eqnarray*}
Let $\tau : \CB(\CH) \to \CB(\CH)$ be defined by $\tau (X) =
\sum_{i=1}^n T_i X T_i^*$ for all $X \in \CB(\CH),$ and let
$\tilde{\tau}^m : M_n (\CB(\CH)) \to M_n (\CB(\CH))$ be defined by
$\tilde{\tau}^m(X)= (\tau^m(X_{ij}))_{n \times n}$ for all $X =
(X_{ij})_{n \times n} \in M_n (\CB(\CH)).$ As $\tau$ is a
completely positive map, $\tilde{\tau}^m$ is also a completely
positive map. \\

So we have $\tilde{\tau}^m(D) \leq I$ and

\begin{eqnarray*}
 \|\tilde{x}_m \|^2 &=& \sum_{r=1}^n \langle \tilde{\tau}^m(D)
 (h_{r1}\ldots h_{rn}) , (h_{r1}\ldots
 h_{rn})\rangle\\
 & \leq & \sum_{r=1}^n \langle (h_{r1}\ldots h_{rn}) , (h_{r1}\ldots
 h_{rn})\rangle\\
 & = & \sum_{r,i}^n \langle h_{ri} , h_{ri}\rangle
  =\sum_{i,r=1}^n \langle T_i^* h_r - \oq_{ir} T_r^* h_i,
 T_i^* h_r - \oq_{ir} T_r^* h_i \rangle\\
 & =&\sum_{i,r=1}^n \{ \langle T_i^*T_ih_r - T_r^*T_ih_i,
 h_r\rangle - \langle T_i^*T_rh_r- T_r^*T_rh_i ,
 h_i\rangle \}\\
& =  &
   \sum_{r=1}^n \langle h_r -\sum_{i=1}^n T_r^*T_ih_i,
 h_r\rangle - \sum_{i=1}^n \langle \sum_{r=1}^n T_i^*T_rh_r-h_i ,
 h_i\rangle \\
 & =&  2 \sum_{r=1}^n \langle h_r -\sum_{i=1}^n T_r^*T_ih_i,
 h_r\rangle
 = 2 \langle D(h_1,\ldots  ,h_n),(h_1,\ldots  ,h_n) \rangle
= 2 \| \tilde{x}_0 \| ^2 =2.
\end{eqnarray*}
As $\langle y, \tilde{x}_0 +\tilde{x}_1\rangle =0$ and $\langle y,
\tilde{x}_{m-1} -\tilde{x}_m \rangle =0 $ for $m + 1 \in
{\mathbb{N}}$,
 we get $\langle y, \tilde{x}_0 +\tilde{x}_{m} \rangle =0$
for $m \in {\mathbb{N}}$. So $1=\langle \tilde {y}_0,
\tilde {y_0}\rangle = \langle \tilde{y}_0,\tilde{x}_0 \rangle = -
\langle \tilde{y}_{m}, \tilde{x}_{m} \rangle$. By Cauchy-Schwarz
inequality, $1 \leq \| \tilde{y}_m \| \| \tilde{x}_m \|$ , which implies
 $\frac{1}{\sqrt{2}} \leq \| \tilde{y}_m \|$ for $m \in
{\mathbb{N}}$. This is a contradiction as $y=0\oplus \oplus
_{m\geq 0}\tilde {y}_m$ is in the  Hilbert space $\tilde {\CH }$.
This proves the particular case.

Using arguements similar to that of Theorem 13 of [BBD], the proof
of the general case (that is when $T_i$ is not necessarily normal)
and the proof of ``$\tilde{\uV}$ is the standard noncommuting
dilation of $\uT$", both follows . \qed

\end{section}

\begin{section}{Distribution of $S_i + S_i^*$ and Related Operator Spaces}
\setcounter{equation}{0}

Let $\CR$ be the von Neumann algebra generated by $G_i=S_i +
S_i^*$ for all $1 \leq i \leq n.$  We are interested in
calculating the moments of $S_i + S_i^*$ with respect to the
vaccum state and inferring about the distribution. The vacuum
expectation is given by $\epsilon (T)=\langle \omega,T \omega
\rangle$ where $T \in \CR.$ So,

$$\epsilon ((S_i + S_i^*)^n)= \langle \omega, (S_i + S_i^*
)^n \omega  \rangle = \left \{ \begin{array}{ccc}
0  & \mbox {~if n is odd} \\
C_n = \frac{1}{n+1} (^n_{\frac{n}{2}}) & \mbox{~otherwise}
\end{array}\right .$$
where $C_n $ the catalan number (refer [Com]). This shows that
$S_i + S_i^*$ has semicircular distribution. Further this vaccum
expectation is not tracial on $\CR$ for $n \geq 2$ as
\begin{eqnarray*}
\epsilon(G_2G_2G_1G_1)&=&\langle \omega, (S_2 + S_2^* )(S_2 +
S_2^*)(S_1 + S_1^* )(S_1 + S_1^*) \omega \rangle\\
&=&\langle \omega, S_2^* S_2^*S_1S_1 + S_2^* S_2S_1^*S_1\omega \rangle = 1 \\
\epsilon(G_2G_1G_1G_2) &=&\langle \omega, (S_2 +
S_2^*)(S_1 + S_1^* )(S_1 + S_1^*)(S_2 + S_2^* ) \omega \rangle\\
&=&\langle \omega, S_2^* S_1^*S_1S_2 + S_2^* S_1S_1^*S_2\omega
\rangle = \frac{1}{2}
\end{eqnarray*}

We would now investigate using arguements of theory of operator
spaces introduced by Effros and Ruan [ER]. Here we follow the
ideas of [BS2] and [HP]. Operator spaces which are Hilbert spaces
are called Hilbertian operator spaces. For some Hilbert space
$\tilde{\CH}$ and $a_i \in B(\tilde{\CH}), 1\leq i \leq n $ define
$$ \|(a_1,\cdots,a_n)\|_{max} = \mbox{max} (\| \sum^n_{i=1} a_ia^*_i\|^{\frac{1}{2}},
\| \sum^n_{i=1} a^*_ia_i\|^{\frac{1}{2}}).$$ Let us denote the
operator space
$$\left \{ \left( \begin{array}{cccc} r_1& 0&\cdots&0\\
. &&& .\\
. &&& .\\
. &&& .\\
r_n & 0 & \cdots& 0
\end{array}\right) \oplus
 \left (\begin{array}{ccc} r_1& \cdots&r_n\\
0 && 0\\
. && .\\
. && .\\
. && .\\
0 & \cdots& 0 \end{array}\right) | r_1,\cdots, r_n \in \C \right\}
\subset M_n \oplus M_n$$ by $E_n.$ Let $\{e_{ij}:1 \leq i,j \leq n
\}$ denote the standard basis of $M_n$ and $\delta_i=e_{i1}\oplus
e_{1i}.$ Then one has
$$ \|\sum^n_{i=1}a_i\otimes \delta_i\|_{B(\tilde{\CH})\otimes M_n}=\|(a_1,\cdots,a_n)\|_{max}.$$
\begin{Theorem}
The operator space generated by $G_i, ~~1 \leq i \leq n$ is
completely isomorphic to $E_n.$
\end{Theorem}
\noindent{\sc Proof:} Its enough to show that for $a_i \in
B(\tilde{\CH}), 1\leq i \leq n$ we have
$$ \|(a_1,\cdots,a_n)\|_{max}\leq \| \sum^n_{i=1}
a_i \otimes G_i\|_{\tilde{\CH}\otimes \Gamma_q(\mathbb{C}^n)} \leq
2 \|(a_1,\cdots,a_n)\|_{max}$$
\begin{eqnarray*}
\| \sum^n_{i=1} a_i \otimes S^*_i\|_{\tilde{\CH}\otimes
\Gamma_q(\mathbb{C}^n)}&=&\| \sum^n_{i=1} (a_i \otimes 1)(1
\otimes
S^*_i)\|_{\tilde{\CH}\otimes \Gamma_q(\mathbb{C}^n)}\\
& \leq & \| \sum^n_{i=1} a_i a^*_i\|^{\frac{1}{2}}_{\tilde{\CH}}
\|\sum^n_{i=1} S_i S^*_i\|^{\frac{1}{2}}_{\Gamma_q(\mathbb{C}^n)}
\leq \| \sum^n_{i=1} a_i
a^*_i\|^{\frac{1}{2}}_{\Gamma_q(\mathbb{C}^n)}
\end{eqnarray*}
Similarly
\begin{eqnarray*}
\|\sum^n_{i=1} a_i \otimes S_i\|_{\tilde{\CH}\otimes
\Gamma_q(\mathbb{C}^n)}&=&\| \sum^n_{i=1} (1 \otimes
S_i)(a_i \otimes 1)\|_{\tilde{\CH}\otimes \Gamma_q(\mathbb{C}^n)}\\
&\leq &  \| \sum^n_{i=1} a^*_i
a_i\|^{\frac{1}{2}}_{\Gamma_q(\mathbb{C}^n)}
\end{eqnarray*}
So
$$  \| \sum^n_{i=1}a_i \otimes G_i\|_{\tilde{\CH}\otimes
\Gamma_q(\mathbb{C}^n)} \leq 2 \|(a_1,\cdots,a_n)\|_{max}$$ Let
$\CS$ denote the set of all states on $B(\tilde{\CH}).$ Now using
the fact that $\epsilon(G_iG_j)=\langle \omega, S^*_i S_j\omega
\rangle= \delta_{ij} $ we get
\begin{eqnarray*}
\| \sum^n_{i=1}a_i \otimes G_i\|^2_{\tilde{\CH}\otimes
\Gamma_q(\mathbb{C}^n)} & \geq & \underbrace{\mbox{sup}}_{\tau \in
\CS} (\tau \otimes \epsilon)[(\sum^n_{i=1}a_i \otimes
G_i)^*\sum^n_{j=1}a_j \otimes G_j]\\
& = & \underbrace{\mbox{sup}}_{\tau \in \CS} \tau (\sum^n_{i=1}
a^*_ia_i) = \| \sum^n_{i=1} a^*_ia_i\|
\end{eqnarray*}

Using similar arguements
$$\| \sum^n_{i=1}a_i \otimes
G_i\|^2_{\tilde{\CH}\otimes \Gamma_q(\mathbb{C}^n)} \geq \|
\sum^n_{i=1} a_ia^*_i\|.$$ \qed

 \vsp \noindent{\bf Acknowledgements:} The author is
supported by a research fellowship from the Indian Statistical
Institute. The author is thankful to B. V. Rajarama Bhat and
Tirthankar Bhattacharyya for many helpful discussions.

\end{section}

\begin{center}
{\bf References}
\end{center}
\begin {itemize}

\item [{[A]}]   J. Agler: The Arveson extension theorem and
coanalytic models, Integral Equations and Operator Theory,
 {\bf 5} (1982), 608-631, {\bf MR 84g:47011}.

\item [{[AP1]}] Arias, A.;  Popescu, G.: Noncommutative interpolation and
Poisson transforms, Israel J. Math., {\bf 115}(2000), 205-234, {\bf
MR 2001i:47021}.

\item [{[Ar1]}] Arveson, W. B.: {\em An Invitation to $C^*$-algebras,\/ }
Graduate Texts in Mathematics, No. 39, Springer-Verlag, New York-Heidelberg
(1976). {\bf MR 58$\#$23621}.

\item [{[Ar2]}]  Arveson, W. B.:   Subalgebras of $C^*$-algebras III, Multivariable
operator theory, Acta Math., {\bf 181}(1998), no. 2, 159-228. {\bf MR 2000e:47013}.

\item [{[At1]}]  Athavale, A. : On the intertwining of joint isometries,
J. Operator Theory, {\bf 23}(1990), 339-350. {\bf MR 91i:47029}.

\item [{[At2]}]  Athavale, A.:  Model theory on the unit ball in $\mathbb{C}^m$,
J. Operator Theory, {\bf 27}(1992), 347-358. {\bf MR 94i:47011}

\item [{[BB]}] Bhat, B. V. Rajarama;   Bhattacharyya, T. : A model theory
for $q$-commuting contractive tuples, J. Operator Theory, {\bf 47}
 (2002), 97-116.

\item [ { [BBD]}] Bhat, B. V. Rajarama;   Bhattacharyya, T.; Dey,
Santanu : Standard noncommuting and commuting dilations of
commuting tuples, to appear in Trans. Amer. Math. Soc.

\item [ {[B]}] Bhattacharyya, T. : Dilation of contractive tuples: a survey; to appear in
{\em Survey of Analysis and Operator Theory }, Proceedings of CMA, Volume 40.

\item [{[BS1]}] Bo\.zejko, M., Speicher, R.: An example of a
generalized Brownian motion, Commun. Math. Phys. {\bf 137} (1991),
no. 3, 519-531 {\bf MR 92m:46096}.

\item [{[BS2]}] Bo\.zejko, M., Speicher, R.: Completely positive maps on Coxeter groups,
 deformed commutation relations, and operator spaces , Math. Ann. {\bf 300} (1994),
  no. 1, 97-120 {\bf MR 95g:46105}.

\item [{[Bu]}] Bunce, J. W.: Models for $n$-tuples of noncommuting
operators, J. Funct. Anal., {\bf 57}(1984), 21-30. {\bf MR
85k:47019}.

 \item [{[Con]}]  Connes, A.:{\em  Noncommutative Geometry,}
Academic Press, (1994) {\bf MR 95j:46063}.

\item [{[Com]}]  Comtet, L.:{\em  Advanced combinatorics,} D. Reidel Publishing Co., (1974)
{\bf MR 57 \#124}.

\item [{[Cu]}] Cuntz, J. : Simple $C^*$-algebras generated by isometries,
Commun. Math. Phys., {\bf 57} (1977),173-185. {\bf MR 57$\# $7189}.

\item [{[Da]}] Davis, C.: Some dilation and representation theorems.
Proceedings of the Second International Symposium in West Africa on
Functional Analysis and its Applications (Kumasi, 1979), 159-182.
{\bf MR 84e:47012}.

\item [{[DKS]}] Davidson K. R.; Kribs, D. W.; Shpigel, M.E.:
Isometric dilations of non-commuting finite rank $n$-tuples,
Canad. J. Math., {\bf 53} (2001) 506-545. {\bf MR 2002f:47010}.

\item [{[ER]}] Effros, E. G.; Ruan, Z. J.:
{\em Operator spaces,\/} London Mathematical Society Monographs.
New Series, {\bf 23},(2000) {\bf MR  2002a:46082}.

\item [{[Fr]}]  Frazho, A. E.: Models for noncommuting operators, J.
Funct. Anal., {\bf 48} (1982), 1-11. {\bf MR  84h:47010}.

\item [{[Ha]}]  Halmos, P. R.: {\em A Hilbert Space Problem Book,\/} Second Edition,
Graduate Texts in Mathematics {\bf 19}, (Springer-Verlag, New
York-Berlin 1982) {\bf MR 84e:47001}.

\item [{[HP]}]  Haagerup, U.; Pisier, G.: Bounded linear operators between
$C *$-algebras, Duke Math. J.,  {\bf 71} (1993), no. 3, 889--925.
{\bf 94k:46112}.

\item [{[JSW]}]  Jorgensen, P. E. T.; Schmitt, L. M.; Werner, R.
F.: $q$-canonical commutation relations and stability of the Cuntz
algebra, Pacific J. Math. {\bf 165} (1994), no. 1, 131-151. {\bf
MR 95g:46116}.

\item [{[M]}]  Majid, S.:  {\em Foundations of Quantum Group Theory,\/} Cambridge
University Press, (1995) {\bf MR 97g:17016 }.

\item [{[Pa]}] Parrott, S.: Unitary dilations for commuting contractions,
Pacific J. Math. {\bf 34} (1970), 481-490. {\bf MR 42$\#$3607}.

\item [{[Po1]}]  Popescu, G.:  Isometric dilations for infinite
sequences of noncommuting operators, Trans. Amer. Math. Soc., {\bf
316}(1989), 523-536. {\bf MR 90c:47006}.

\item [{[Po2]}]
 Popescu, G.: Characteristic functions for infinite
sequences of noncommuting operators, J. Operator Theory, {\bf 22} (1989),
51-71. {\bf MR 91m:47012}.

\item [{[Po3]}] Popescu, G.: Poisson transforms on some $C\sp *$-algebras
generated by isometries. J. Funct. Anal. 161 (1999), 27-61. {\bf MR 2000m:46117}.

\item [{[Po4]}]  Popescu, G.:  Curvature invariant for Hilbert modules
over free semigroup algebras, Advances in  Mathematics,
{\bf 158} (2001), 264-309. {\bf MR  2002b:46097}.

\item [{[Pr]}]  Prugovecki, E.: {\em Quantum Mechanics in Hilbert Space, \/}
Academic Press, Second Edition(1981) {\bf MR 84k:81005 }.

\item [{[Pu]}]  Putnam, C. R.: {\em Commutation properties of Hilbert space
Operators and Related Topics,\/} Springer-Verlag (1967). {\bf MR 36$\#$707}.


\end{itemize}

 \noindent {\sc Santanu Dey} \\ Indian Statistical Institute,
R. V. College Post, Bangalore 560059, India.\\
 e-mail: {\sl santanu@isibang.ac.in }\\

\end{document}